%% file: manuscript.tex
\documentclass[final,hidelinks,onefignum,onetabnum]{siamart250211}

\input{ex_shared}

\ifpdf
\hypersetup{
	pdftitle={Optimal Design of Broadband Absorbers with Multiple Plasmonic Nanoparticles via Reduced Basis Method},
	pdfauthor={Y. Gao, H. Zhang and K. Zhang}
}
\fi

\begin{document}

\maketitle

% REQUIRED
\begin{abstract}
	In this paper, we propose a computational framework for the optimal design of broadband absorbing materials composed of plasmonic nanoparticle arrays. This design problem poses several key challenges: (1) the complex multi-particle interactions and high-curvature geometries; (2) the requirement to achieve broadband frequency responses, including resonant regimes; (3) the complexity of shape derivative calculations; and (4) the non-convexity of the optimization landscape. To systematically address these challenges, we employ three sequential strategies. First, we introduce a parameterized integral equation formulation that circumvents traditional shape derivative computations. Second, we develop a shape-adaptive reduced basis method (RBM) that utilizes the eigenfunctions of the Neumann–Poincar\'{e} operator for forward problems and their adjoint counterparts for adjoint problems, thereby addressing singularities and accelerating computations. Third, we propose a physics-informed initialization strategy that estimates nanoparticle configurations under weak coupling assumptions, thereby improving the performance of gradient-based optimization algorithms. The method's computational advantages are demonstrated through numerical experiments, which show accurate and efficient designs across various geometric configurations. Furthermore, the framework is flexible and extensible to other material systems and boundary conditions.
\end{abstract}

% REQUIRED

% REQUIRED

\section{Introduction}

 Controlling light absorption at the nanoscale poses significant scientific and engineering challenges. In recent years, noble metal nanoparticles have emerged as promising candidates due to their exceptional plasmonic properties and their ability to efficiently confine light. Consequently, the optimal design of broadband absorbers incorporating multiple plasmonic nanoparticles has become an important topic in nanophotonics. This growing interest is largely motivated by applications in solar energy harvesting, thermal emission control, and the enhancement of various optoelectronic devices (see~\cite{estakhri2013physics,grigoriev2015optimizing,riley2017near,wang2023review}).

The absorption mechanism in metals originates from the positive imaginary part of their permittivity, while the negative real part enables the excitation of localized surface plasmon resonances. At these resonant frequencies, both the electromagnetic field and the absorptance are significantly enhanced. In the subwavelength regime, plasmonic phenomena can be rigorously analyzed using the quasi-static approximation~\cite{ammari2016surface,ammari2017mathematical,ammari2018heat,ammari2016mathematical,yu2018plasmonic}. This framework reduces the resonance analysis to the spectral study of the Neumann--Poincar\'{e} (NP) operator.
As demonstrated in~\cite{ammari2016surface,ammari2017mathematical}, the resonant frequencies are dependent on the geometry of the nanoparticles. For example, a circular nanoparticle exhibits a single narrow absorptance peak, whereas an elliptical nanoparticle supports two distinct resonances due to its anisotropic shape. Consequently, achieving broadband absorptance requires an ensemble of nanoparticles, each contributing a narrow spectral response at resonant frequencies, thereby collectively spanning a broad bandwidth. In this study, we investigate the optimal design of broadband absorbers by multiple plasmonic nanoparticles. 

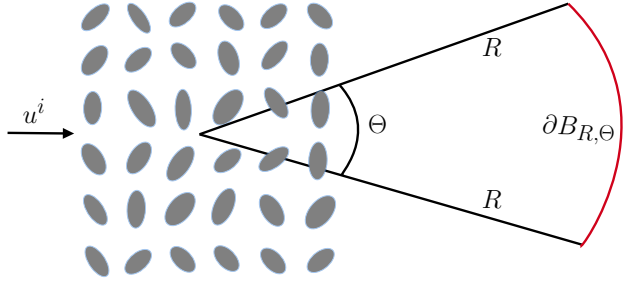
\begin{figure}[htbp]
	\begin{center}
		\scalebox{0.40}{
		\tikzset{every picture/.style={line width=0.75pt}} %set default line width to 0.75pt        

\begin{tikzpicture}[x=0.75pt,y=0.75pt,yscale=-1,xscale=1]
	%uncomment if require: \path (0,395); %set diagram left start at 0, and has height of 395
	
	%Shape: Ellipse [id:dp7817467131910729] 
	\draw  [color={rgb, 255:red, 74; green, 144; blue, 226 }  ,draw opacity=0.5 ][fill={rgb, 255:red, 128; green, 128; blue, 128 }  ,fill opacity=1 ] (175.08,225.61) .. controls (182.69,216.09) and (192.27,212.9) .. (196.48,218.5) .. controls (200.7,224.11) and (197.95,236.37) .. (190.34,245.9) .. controls (182.74,255.43) and (173.15,258.61) .. (168.94,253.01) .. controls (164.73,247.41) and (167.48,235.14) .. (175.08,225.61) -- cycle ;
	%Shape: Ellipse [id:dp7805628975234729] 
	\draw  [color={rgb, 255:red, 74; green, 144; blue, 226 }  ,draw opacity=0.5 ][fill={rgb, 255:red, 128; green, 128; blue, 128 }  ,fill opacity=1 ] (170.21,152.2) .. controls (174.55,146.77) and (185.21,151.85) .. (194.01,163.55) .. controls (202.81,175.25) and (206.43,189.15) .. (202.09,194.58) .. controls (197.75,200.02) and (187.1,194.94) .. (178.29,183.24) .. controls (169.49,171.53) and (165.87,157.64) .. (170.21,152.2) -- cycle ;
	%Shape: Ellipse [id:dp7097967538417123] 
	\draw  [color={rgb, 255:red, 74; green, 144; blue, 226 }  ,draw opacity=0.5 ][fill={rgb, 255:red, 128; green, 128; blue, 128 }  ,fill opacity=1 ] (286.93,159.38) .. controls (296.57,149.14) and (307.99,146.41) .. (312.46,153.28) .. controls (316.93,160.16) and (312.74,174.03) .. (303.11,184.27) .. controls (293.48,194.51) and (282.05,197.24) .. (277.58,190.37) .. controls (273.12,183.5) and (277.3,169.62) .. (286.93,159.38) -- cycle ;
	%Shape: Ellipse [id:dp9448633806413029] 
	\draw  [color={rgb, 255:red, 74; green, 144; blue, 226 }  ,draw opacity=0.5 ][fill={rgb, 255:red, 128; green, 128; blue, 128 }  ,fill opacity=1 ] (123.87,151.97) .. controls (130.14,151.23) and (135.4,159.98) .. (135.6,171.5) .. controls (135.8,183.02) and (130.88,192.96) .. (124.6,193.69) .. controls (118.33,194.43) and (113.08,185.69) .. (112.87,174.16) .. controls (112.67,162.64) and (117.59,152.7) .. (123.87,151.97) -- cycle ;
	%Shape: Ellipse [id:dp6977257685678746] 
	\draw  [color={rgb, 255:red, 74; green, 144; blue, 226 }  ,draw opacity=0.5 ][fill={rgb, 255:red, 128; green, 128; blue, 128 }  ,fill opacity=1 ] (117.16,246.85) .. controls (109.86,237.28) and (107.74,225.29) .. (112.42,220.08) .. controls (117.1,214.87) and (126.82,218.41) .. (134.12,227.99) .. controls (141.42,237.57) and (143.54,249.55) .. (138.86,254.76) .. controls (134.18,259.97) and (124.46,256.43) .. (117.16,246.85) -- cycle ;
	%Shape: Ellipse [id:dp8432279325231655] 
	\draw  [color={rgb, 255:red, 74; green, 144; blue, 226 }  ,draw opacity=0.5 ][fill={rgb, 255:red, 128; green, 128; blue, 128 }  ,fill opacity=1 ] (223.83,94.39) .. controls (226.7,88.55) and (235.76,89.34) .. (244.08,96.16) .. controls (252.39,102.97) and (256.8,113.23) .. (253.94,119.06) .. controls (251.07,124.9) and (242.01,124.11) .. (233.7,117.29) .. controls (225.38,110.48) and (220.97,100.22) .. (223.83,94.39) -- cycle ;
	%Shape: Ellipse [id:dp16764489545513173] 
	\draw  [color={rgb, 255:red, 74; green, 144; blue, 226 }  ,draw opacity=0.5 ][fill={rgb, 255:red, 128; green, 128; blue, 128 }  ,fill opacity=1 ] (120.26,101.48) .. controls (128.81,93.16) and (138.66,91.42) .. (142.27,97.6) .. controls (145.88,103.79) and (141.87,115.55) .. (133.32,123.88) .. controls (124.77,132.2) and (114.91,133.94) .. (111.3,127.76) .. controls (107.7,121.57) and (111.71,109.81) .. (120.26,101.48) -- cycle ;
	%Shape: Ellipse [id:dp665170948719304] 
	\draw  [color={rgb, 255:red, 74; green, 144; blue, 226 }  ,draw opacity=0.5 ][fill={rgb, 255:red, 128; green, 128; blue, 128 }  ,fill opacity=1 ] (237.5,149.23) .. controls (243.15,148.82) and (248.2,159.95) .. (248.77,174.1) .. controls (249.35,188.25) and (245.24,200.06) .. (239.59,200.47) .. controls (233.94,200.88) and (228.89,189.74) .. (228.32,175.59) .. controls (227.74,161.44) and (231.85,149.64) .. (237.5,149.23) -- cycle ;
	%Straight Lines [id:da03789014046385364] 
	\draw [line width=2.25]    (17,204.11) -- (96,204.89) ;
	\draw [shift={(101,204.94)}, rotate = 180.57] [fill={rgb, 255:red, 0; green, 0; blue, 0 }  ][line width=0.08]  [draw opacity=0] (14.29,-6.86) -- (0,0) -- (14.29,6.86) -- cycle    ;
	%Shape: Ellipse [id:dp18866677010312127] 
	\draw  [color={rgb, 255:red, 74; green, 144; blue, 226 }  ,draw opacity=0.5 ][fill={rgb, 255:red, 128; green, 128; blue, 128 }  ,fill opacity=1 ] (164.85,119.86) .. controls (165.05,112.43) and (172.73,102.83) .. (182.02,98.41) .. controls (191.3,93.98) and (198.66,96.42) .. (198.46,103.84) .. controls (198.26,111.26) and (190.58,120.87) .. (181.3,125.29) .. controls (172.02,129.72) and (164.65,127.28) .. (164.85,119.86) -- cycle ;
	%Shape: Ellipse [id:dp9108829278572967] 
	\draw  [color={rgb, 255:red, 74; green, 144; blue, 226 }  ,draw opacity=0.5 ][fill={rgb, 255:red, 128; green, 128; blue, 128 }  ,fill opacity=1 ] (178.73,324.27) .. controls (172.31,324.51) and (167.35,314.04) .. (167.67,300.9) .. controls (167.98,287.76) and (173.44,276.91) .. (179.87,276.67) .. controls (186.3,276.43) and (191.25,286.89) .. (190.94,300.04) .. controls (190.62,313.18) and (185.16,324.03) .. (178.73,324.27) -- cycle ;
	%Shape: Ellipse [id:dp5490228337862861] 
	\draw  [color={rgb, 255:red, 74; green, 144; blue, 226 }  ,draw opacity=0.5 ][fill={rgb, 255:red, 128; green, 128; blue, 128 }  ,fill opacity=1 ] (134.84,293.28) .. controls (142.17,302.97) and (144.94,314.4) .. (141.02,318.81) .. controls (137.1,323.22) and (127.98,318.93) .. (120.65,309.24) .. controls (113.32,299.55) and (110.56,288.12) .. (114.48,283.71) .. controls (118.4,279.3) and (127.52,283.59) .. (134.84,293.28) -- cycle ;
	%Shape: Ellipse [id:dp30846502414533217] 
	\draw  [color={rgb, 255:red, 74; green, 144; blue, 226 }  ,draw opacity=0.5 ][fill={rgb, 255:red, 128; green, 128; blue, 128 }  ,fill opacity=1 ] (250.25,223.26) .. controls (254.82,228.53) and (251.67,241.61) .. (243.22,252.47) .. controls (234.78,263.33) and (224.23,267.85) .. (219.66,262.58) .. controls (215.1,257.3) and (218.24,244.22) .. (226.69,233.36) .. controls (235.14,222.51) and (245.69,217.98) .. (250.25,223.26) -- cycle ;
	%Shape: Ellipse [id:dp7528388513042032] 
	\draw  [color={rgb, 255:red, 74; green, 144; blue, 226 }  ,draw opacity=0.5 ][fill={rgb, 255:red, 128; green, 128; blue, 128 }  ,fill opacity=1 ] (310.58,227.62) .. controls (313.51,233.12) and (308.59,242.74) .. (299.6,249.09) .. controls (290.61,255.45) and (280.96,256.13) .. (278.04,250.63) .. controls (275.11,245.13) and (280.03,235.51) .. (289.02,229.16) .. controls (298.01,222.8) and (307.66,222.11) .. (310.58,227.62) -- cycle ;
	%Shape: Ellipse [id:dp4286506045351679] 
	\draw  [color={rgb, 255:red, 74; green, 144; blue, 226 }  ,draw opacity=0.5 ][fill={rgb, 255:red, 128; green, 128; blue, 128 }  ,fill opacity=1 ] (251.87,283.26) .. controls (256.3,289.98) and (252.11,303.05) .. (242.53,312.45) .. controls (232.94,321.86) and (221.57,324.04) .. (217.14,317.32) .. controls (212.71,310.6) and (216.89,297.53) .. (226.48,288.12) .. controls (236.07,278.72) and (247.44,276.54) .. (251.87,283.26) -- cycle ;
	%Shape: Ellipse [id:dp32337057581429907] 
	\draw  [color={rgb, 255:red, 74; green, 144; blue, 226 }  ,draw opacity=0.5 ][fill={rgb, 255:red, 128; green, 128; blue, 128 }  ,fill opacity=1 ] (300.52,275.4) .. controls (306.09,279.41) and (305.57,292.49) .. (299.36,304.62) .. controls (293.14,316.75) and (283.6,323.32) .. (278.04,319.31) .. controls (272.47,315.3) and (272.99,302.22) .. (279.2,290.09) .. controls (285.41,277.97) and (294.96,271.39) .. (300.52,275.4) -- cycle ;
	%Shape: Ellipse [id:dp6496084449913768] 
	\draw  [color={rgb, 255:red, 74; green, 144; blue, 226 }  ,draw opacity=0.5 ][fill={rgb, 255:red, 128; green, 128; blue, 128 }  ,fill opacity=1 ] (286.72,89.67) .. controls (292.36,85.46) and (300.97,91.11) .. (305.93,102.29) .. controls (310.9,113.47) and (310.35,125.95) .. (304.71,130.17) .. controls (299.07,134.38) and (290.46,128.73) .. (285.5,117.54) .. controls (280.53,106.36) and (281.08,93.88) .. (286.72,89.67) -- cycle ;
	%Straight Lines [id:da8134946665298426] 
	\draw [line width=2.25]    (741,345.11) -- (259.04,205.94) ;
	%Straight Lines [id:da02735840097660236] 
	\draw [line width=2.25]    (724,41.11) -- (259.04,205.94) ;
	%Shape: Ellipse [id:dp5199812230275542] 
	\draw  [color={rgb, 255:red, 74; green, 144; blue, 226 }  ,draw opacity=0.5 ][fill={rgb, 255:red, 128; green, 128; blue, 128 }  ,fill opacity=1 ] (358.84,292.28) .. controls (366.17,301.97) and (368.26,314.16) .. (363.52,319.5) .. controls (358.77,324.84) and (348.98,321.31) .. (341.65,311.61) .. controls (334.32,301.92) and (332.23,289.73) .. (336.98,284.4) .. controls (341.73,279.06) and (351.52,282.59) .. (358.84,292.28) -- cycle ;
	%Shape: Ellipse [id:dp35202907242287473] 
	\draw  [color={rgb, 255:red, 74; green, 144; blue, 226 }  ,draw opacity=0.5 ][fill={rgb, 255:red, 128; green, 128; blue, 128 }  ,fill opacity=1 ] (344.45,230.54) .. controls (353.8,222.77) and (365.09,219.49) .. (369.67,223.2) .. controls (374.25,226.91) and (370.39,236.22) .. (361.04,243.98) .. controls (351.7,251.75) and (340.41,255.03) .. (335.82,251.32) .. controls (331.24,247.61) and (335.1,238.3) .. (344.45,230.54) -- cycle ;
	%Shape: Ellipse [id:dp7531615376551197] 
	\draw  [color={rgb, 255:red, 74; green, 144; blue, 226 }  ,draw opacity=0.5 ][fill={rgb, 255:red, 128; green, 128; blue, 128 }  ,fill opacity=1 ] (363.3,111.03) .. controls (355.85,120.63) and (345.53,126.28) .. (340.26,123.64) .. controls (334.98,121.01) and (336.74,111.09) .. (344.19,101.49) .. controls (351.64,91.89) and (361.96,86.24) .. (367.24,88.88) .. controls (372.51,91.51) and (370.75,101.43) .. (363.3,111.03) -- cycle ;
	%Shape: Ellipse [id:dp3496887763139722] 
	\draw  [color={rgb, 255:red, 74; green, 144; blue, 226 }  ,draw opacity=0.5 ][fill={rgb, 255:red, 128; green, 128; blue, 128 }  ,fill opacity=1 ] (362.84,158.28) .. controls (370.17,167.97) and (372.94,179.4) .. (369.02,183.81) .. controls (365.1,188.22) and (355.98,183.93) .. (348.65,174.24) .. controls (341.32,164.55) and (338.56,153.12) .. (342.48,148.71) .. controls (346.4,144.3) and (355.52,148.59) .. (362.84,158.28) -- cycle ;
	%Shape: Ellipse [id:dp16889733499433168] 
	\draw  [color={rgb, 255:red, 74; green, 144; blue, 226 }  ,draw opacity=0.5 ][fill={rgb, 255:red, 128; green, 128; blue, 128 }  ,fill opacity=1 ] (136.84,358.28) .. controls (144.17,367.97) and (146.94,379.4) .. (143.02,383.81) .. controls (139.1,388.22) and (129.98,383.93) .. (122.65,374.24) .. controls (115.32,364.55) and (112.56,353.12) .. (116.48,348.71) .. controls (120.4,344.3) and (129.52,348.59) .. (136.84,358.28) -- cycle ;
	%Shape: Ellipse [id:dp1358295975631605] 
	\draw  [color={rgb, 255:red, 74; green, 144; blue, 226 }  ,draw opacity=0.5 ][fill={rgb, 255:red, 128; green, 128; blue, 128 }  ,fill opacity=1 ] (118.72,52.55) .. controls (127.11,43.75) and (137.94,39.17) .. (142.92,42.33) .. controls (147.91,45.49) and (145.15,55.18) .. (136.77,63.98) .. controls (128.39,72.77) and (117.55,77.35) .. (112.57,74.19) .. controls (107.59,71.04) and (110.34,61.34) .. (118.72,52.55) -- cycle ;
	%Shape: Ellipse [id:dp04955051780406272] 
	\draw  [color={rgb, 255:red, 74; green, 144; blue, 226 }  ,draw opacity=0.5 ][fill={rgb, 255:red, 128; green, 128; blue, 128 }  ,fill opacity=1 ] (188.68,52.72) .. controls (196.01,62.41) and (199.25,73.3) .. (195.93,77.03) .. controls (192.61,80.76) and (183.98,75.93) .. (176.65,66.24) .. controls (169.32,56.55) and (166.07,45.66) .. (169.39,41.93) .. controls (172.72,38.19) and (181.35,43.02) .. (188.68,52.72) -- cycle ;
	%Shape: Ellipse [id:dp561740039944477] 
	\draw  [color={rgb, 255:red, 74; green, 144; blue, 226 }  ,draw opacity=0.5 ][fill={rgb, 255:red, 128; green, 128; blue, 128 }  ,fill opacity=1 ] (220.85,62.86) .. controls (221.05,55.43) and (228.73,45.83) .. (238.02,41.41) .. controls (247.3,36.98) and (254.66,39.42) .. (254.46,46.84) .. controls (254.26,54.26) and (246.58,63.87) .. (237.3,68.29) .. controls (228.02,72.72) and (220.65,70.28) .. (220.85,62.86) -- cycle ;
	%Shape: Ellipse [id:dp05040463790357186] 
	\draw  [color={rgb, 255:red, 74; green, 144; blue, 226 }  ,draw opacity=0.5 ][fill={rgb, 255:red, 128; green, 128; blue, 128 }  ,fill opacity=1 ] (164.85,374.86) .. controls (165.05,367.43) and (172.73,357.83) .. (182.02,353.41) .. controls (191.3,348.98) and (198.66,351.42) .. (198.46,358.84) .. controls (198.26,366.26) and (190.58,375.87) .. (181.3,380.29) .. controls (172.02,384.72) and (164.65,382.28) .. (164.85,374.86) -- cycle ;
	%Shape: Ellipse [id:dp5114974191202157] 
	\draw  [color={rgb, 255:red, 74; green, 144; blue, 226 }  ,draw opacity=0.5 ][fill={rgb, 255:red, 128; green, 128; blue, 128 }  ,fill opacity=1 ] (335.85,64.86) .. controls (336.05,57.43) and (343.73,47.83) .. (353.02,43.41) .. controls (362.3,38.98) and (369.66,41.42) .. (369.46,48.84) .. controls (369.26,56.26) and (361.58,65.87) .. (352.3,70.29) .. controls (343.02,74.72) and (335.65,72.28) .. (335.85,64.86) -- cycle ;
	%Shape: Ellipse [id:dp37488347771952135] 
	\draw  [color={rgb, 255:red, 74; green, 144; blue, 226 }  ,draw opacity=0.5 ][fill={rgb, 255:red, 128; green, 128; blue, 128 }  ,fill opacity=1 ] (222.83,351.39) .. controls (225.7,345.55) and (234.76,346.34) .. (243.08,353.16) .. controls (251.39,359.97) and (255.8,370.23) .. (252.94,376.06) .. controls (250.07,381.9) and (241.01,381.11) .. (232.7,374.29) .. controls (224.38,367.48) and (219.97,357.22) .. (222.83,351.39) -- cycle ;
	%Shape: Ellipse [id:dp9819940723136769] 
	\draw  [color={rgb, 255:red, 74; green, 144; blue, 226 }  ,draw opacity=0.5 ][fill={rgb, 255:red, 128; green, 128; blue, 128 }  ,fill opacity=1 ] (278.04,350.39) .. controls (280.9,344.55) and (289.96,345.34) .. (298.28,352.16) .. controls (306.59,358.97) and (311.01,369.23) .. (308.14,375.06) .. controls (305.27,380.9) and (296.21,380.11) .. (287.9,373.29) .. controls (279.58,366.48) and (275.17,356.22) .. (278.04,350.39) -- cycle ;
	%Shape: Ellipse [id:dp5928899590636589] 
	\draw  [color={rgb, 255:red, 74; green, 144; blue, 226 }  ,draw opacity=0.5 ][fill={rgb, 255:red, 128; green, 128; blue, 128 }  ,fill opacity=1 ] (284.04,46.39) .. controls (286.9,40.55) and (295.96,41.34) .. (304.28,48.16) .. controls (312.59,54.97) and (317.01,65.23) .. (314.14,71.06) .. controls (311.27,76.9) and (302.21,76.11) .. (293.9,69.29) .. controls (285.58,62.48) and (281.17,52.22) .. (284.04,46.39) -- cycle ;
	%Shape: Ellipse [id:dp9040161896857479] 
	\draw  [color={rgb, 255:red, 74; green, 144; blue, 226 }  ,draw opacity=0.5 ][fill={rgb, 255:red, 128; green, 128; blue, 128 }  ,fill opacity=1 ] (336.83,353.39) .. controls (339.7,347.55) and (348.76,348.34) .. (357.08,355.16) .. controls (365.39,361.97) and (369.8,372.23) .. (366.94,378.06) .. controls (364.07,383.9) and (355.01,383.11) .. (346.7,376.29) .. controls (338.38,369.48) and (333.97,359.22) .. (336.83,353.39) -- cycle ;
	%Shape: Ellipse [id:dp4451898367956122] 
	\draw  [color={rgb, 255:red, 74; green, 144; blue, 226 }  ,draw opacity=0.5 ][fill={rgb, 255:red, 128; green, 128; blue, 128 }  ,fill opacity=1 ] (395.83,44.39) .. controls (398.7,38.55) and (407.76,39.34) .. (416.08,46.16) .. controls (424.39,52.97) and (428.8,63.23) .. (425.94,69.06) .. controls (423.07,74.9) and (414.01,74.11) .. (405.7,67.29) .. controls (397.38,60.48) and (392.97,50.22) .. (395.83,44.39) -- cycle ;
	%Shape: Ellipse [id:dp564782340847837] 
	\draw  [color={rgb, 255:red, 74; green, 144; blue, 226 }  ,draw opacity=0.5 ][fill={rgb, 255:red, 128; green, 128; blue, 128 }  ,fill opacity=1 ] (409.87,90.97) .. controls (416.14,90.23) and (421.4,98.98) .. (421.6,110.5) .. controls (421.8,122.02) and (416.88,131.96) .. (410.6,132.69) .. controls (404.33,133.43) and (399.08,124.69) .. (398.87,113.16) .. controls (398.67,101.64) and (403.59,91.7) .. (409.87,90.97) -- cycle ;
	%Shape: Ellipse [id:dp7044419845093273] 
	\draw  [color={rgb, 255:red, 74; green, 144; blue, 226 }  ,draw opacity=0.5 ][fill={rgb, 255:red, 128; green, 128; blue, 128 }  ,fill opacity=1 ] (410.73,198.27) .. controls (404.31,198.51) and (399.35,188.04) .. (399.67,174.9) .. controls (399.98,161.76) and (405.44,150.91) .. (411.87,150.67) .. controls (418.3,150.43) and (423.25,160.89) .. (422.94,174.04) .. controls (422.62,187.18) and (417.16,198.03) .. (410.73,198.27) -- cycle ;
	%Shape: Ellipse [id:dp8402788334497089] 
	\draw  [color={rgb, 255:red, 74; green, 144; blue, 226 }  ,draw opacity=0.5 ][fill={rgb, 255:red, 128; green, 128; blue, 128 }  ,fill opacity=1 ] (407.73,263.27) .. controls (401.31,263.51) and (396.35,253.04) .. (396.67,239.9) .. controls (396.98,226.76) and (402.44,215.91) .. (408.87,215.67) .. controls (415.3,215.43) and (420.25,225.89) .. (419.94,239.04) .. controls (419.62,252.18) and (414.16,263.03) .. (407.73,263.27) -- cycle ;
	%Shape: Ellipse [id:dp7307656611854811] 
	\draw  [color={rgb, 255:red, 74; green, 144; blue, 226 }  ,draw opacity=0.5 ][fill={rgb, 255:red, 128; green, 128; blue, 128 }  ,fill opacity=1 ] (428.87,287.26) .. controls (433.3,293.98) and (429.11,307.05) .. (419.53,316.45) .. controls (409.94,325.86) and (398.57,328.04) .. (394.14,321.32) .. controls (389.71,314.6) and (393.89,301.53) .. (403.48,292.12) .. controls (413.07,282.72) and (424.44,280.54) .. (428.87,287.26) -- cycle ;
	%Shape: Ellipse [id:dp6305243683512922] 
	\draw  [color={rgb, 255:red, 74; green, 144; blue, 226 }  ,draw opacity=0.5 ][fill={rgb, 255:red, 128; green, 128; blue, 128 }  ,fill opacity=1 ] (428.08,353.62) .. controls (431,359.12) and (426.08,368.74) .. (417.09,375.09) .. controls (408.11,381.45) and (398.45,382.13) .. (395.53,376.63) .. controls (392.61,371.13) and (397.53,361.51) .. (406.51,355.16) .. controls (415.5,348.8) and (425.16,348.11) .. (428.08,353.62) -- cycle ;
	\draw [color={rgb, 255:red, 208; green, 2; blue, 27 }  ,draw opacity=1 ][line width=2.25]    (724,41.11) .. controls (869,183.11) and (732.5,362) .. (741,345.11) ;
	%Curve Lines [id:da5932285500685144] 
	\draw [line width=2.25]    (436,144.11) .. controls (490.5,207) and (429.5,263) .. (440,256.11) ;
	
	% Text Node
	\draw (34.95,160) node [anchor=north west][inner sep=0.75pt]  [font=\Huge]  {$u^{i}$};
	% Text Node
	\draw (690,180) node [anchor=north west][inner sep=0.75pt]  [font=\Huge,rotate=-1.31]  {$\partial B_{R,\Theta }$};
	% Text Node
	\draw (613.5,275) node [anchor=north west][inner sep=0.75pt]  [font=\Huge]  {$R$};
	% Text Node
	\draw (613.5,82.9) node [anchor=north west][inner sep=0.75pt]  [font=\Huge]  {$R$};
	% Text Node
	\draw (470,180) node [anchor=north west][inner sep=0.75pt]  [font=\Huge]  {$\Theta $};

\end{tikzpicture}
	}
	\end{center}
	\caption{Illustration of multiple plasmonic nanoparticles $D= \bigcup_{m=1}^M D_m$ and a partial sphere $\partial B_{R,\Theta}$ (see \eqref{eq:partial_sphere}) for energy measurement.}
	\label{fig:multiple_nanoparticles}
\end{figure}

\subsection{Background}

We investigate electromagnetic scattering in an inhomogeneous medium consisting of plasmonic nanoparticles, as illustrated in \cref{fig:multiple_nanoparticles}. Specifically, the nanoparticle region $D$ comprises well-separated particles whose permittivity $\varepsilon_c(\lambda)$ and permeability $\mu_c(\lambda)$ depend on the wavelength $\lambda$. In contrast, the surrounding homogeneous medium $\mathbb{R}^2 \setminus \overline{D}$ is characterized by constant permittivity $\varepsilon_m$ and permeability $\mu_m$, both independent of $\lambda$. The material parameters are given by
\begin{align*}
    \varepsilon_c(\lambda)  =\varepsilon_0\varepsilon_{r,c}(\lambda),\quad \mu_c(\lambda) =\mu_0\mu_{r,c}(\lambda), \quad\varepsilon_m =\varepsilon_0\varepsilon_{r,m},\quad \mu_m =\mu_0\mu_{r,m},
\end{align*}
where $\varepsilon_0$ and $\mu_0$ denote the permittivity and permeability of vacuum, respectively, and $\varepsilon_{r,\cdot}$ and $\mu_{r,\cdot}$ represent the corresponding relative permittivities and permeabilities. Consequently, the spatial distributions of the material parameters are described by
\begin{align*}
	\varepsilon(\lambda, x) = \varepsilon_c(\lambda) \chi_D(x) + \varepsilon_m \chi_{\mathbb{R}^2 \setminus \overline{D}}(x), \quad
	\mu(\lambda, x) = \mu_c(\lambda) \chi_D(x) + \mu_m \chi_{\mathbb{R}^2 \setminus \overline{D}}(x),
\end{align*}
where $\chi$ is the characteristic function. The wavenumber in vacuum is defined as $k_0(\lambda) = \omega(\lambda) \sqrt{\varepsilon_0 \mu_0} = \omega(\lambda) / c_0 = 2\pi/\lambda$, with $\omega$ denoting the angular frequency and $c_0$ the speed of light in vacuum. Accordingly, the wavenumbers in the homogeneous medium and the nanoparticle regions are given by $k_m(\lambda) = k_0(\lambda) \sqrt{\varepsilon_{r,m} \mu_{r,m}}$ and $k_c(\lambda) = k_0(\lambda) \sqrt{\varepsilon_{r,c}(\lambda) \mu_{r,c}(\lambda)}$, respectively. We suppress $\lambda$-dependence in $k_m$, $k_c$, $\varepsilon$, and $\mu$ hereafter for the ease of notation.

For simplicity,  we restrict our discussion to transverse magnetic (TM) polarization. Let the magnetic field $\mathbf{H}(x_1, x_2) = (0, 0, u(x_1, x_2))^{\top}$, where $u$ is a scalar function. The corresponding multiple scattering problem is governed by the following system
\begin{equation}\label{eq:original_scattering_problem}
	\begin{cases}
		\displaystyle
		\nabla \cdot \left( \frac{1}{\varepsilon} \nabla u \right) + \omega^2 \mu u = 0 
		& \text{in } \mathbb{R}^2 \setminus \partial D, \\[2mm]
		u_{+} - u_{-} = 0 
		& \text{on } \partial D, \\[2mm]
		\displaystyle
		\left. \frac{1}{\varepsilon_m} \frac{\partial u}{\partial \nu} \right|_{+} 
		- \left. \frac{1}{\varepsilon_c} \frac{\partial u}{\partial \nu} \right|_{-} = 0 
		& \text{on } \partial D.
	\end{cases}
\end{equation}
The incident field $u^i(x) = e^{\mathrm{i} k_m x \cdot d}$ propagates along $d = (\cos \theta_0, \sin \theta_0)$ with incident angle $\theta_0$, while the scattered field $u^s = u - u^i$ satisfies the Sommerfeld radiation condition \cite{CK2019Inverse}. In general, the optimal design problem can be formulated as follows: \emph{Given an incident wave $u^i$, the objective is to design an array of nanoparticles arranged around a reference point such that the absorptance $A(\lambda)$ (see~\eqref{eq:absorptance}), measured in the forward direction within a specified angular range $\Theta\subseteq [0,2\pi)$, matches a specified target value $A^{\mathrm{tar}}(\lambda)$ over a broad wavelength range $\lambda\in \Lambda=[\lambda_{\min}, \lambda_{\max}]$.}

This design problem poses substantial theoretical and computational challenges. From a theoretical perspective, fundamental limits for passive linear systems constrain the relationship between absorptance efficiency, bandwidth, and particle geometry (see~\cite{bernland2011sum,cassier2017bounds,rozanov2002ultimate,rozanov2010limitations}, and references therein). Beyond these physical bounds, the broadband optimization of multi-particle systems introduces additional computational challenges, primarily due to:

\begin{itemize}[leftmargin=*]
	\item \textbf{Multiple scattering:} Interactions among particles lead to large-scale systems, with computational complexity growing rapidly as the number of particles $M$ increases. Additionally, the presence of high-curvature features in certain nanoparticles demands rigorous numerical treatment.
	
	\item \textbf{Broadband computation:} Solving multiple scattering over a broad spectral range $[\lambda_{\min}, \lambda_{\max}]$ is computationally expensive, as it requires evaluations at numerous discrete wavelengths $\lambda$. Geometries near plasmonic resonance can induce ill-conditioned system matrices, necessitating specialized solvers to ensure numerical stability and accuracy.
	
	\item \textbf{Shape derivative:} The problem involves multi-domain shape optimization, where shape derivatives must account for both intra-particle and inter-particle interactions, unlike conventional single-domain scattering problems.
	
	\item \textbf{Non-convexity:} The design landscape is highly non-convex, making the identification of global optima challenging. A physics-informed initial guess, guided by multiple scattering theory, is therefore crucial for achieving convergence.
	
	\item \textbf{High dimensionality:} The optimal design problem is inherently high-dimensional, due to the large number of particles in the discretization of both the scattering and adjoint problems, as well as the resulting number of optimization variables.
\end{itemize}

Over the years, numerous numerical methods have been developed for electromagnetic scattering problems, including finite difference methods (FDM,~\cite{kunz1993finite,taflove2005computational}) and finite element methods (FEM,~\cite{jin2015finite,monk2003finite}), and etc. However, the boundary element method (BEM) offers several advantages over these techniques in some cases. First, it represents solutions solely in terms of surface densities, significantly reducing the number of unknowns. Second, it inherently satisfies the radiation condition, thereby eliminating errors associated with artificial boundary conditions~\cite{colton2013integral,kanwal2013linear}. Third, multiple-scattering scenarios introduce additional complexity due to interactions among particles. As inter-scatterer separation increases, FDM and FEM require larger computational domains, creating two computational burdens: (1) the resulting linear systems become increasingly expensive to solve, and (2) material interfaces typically demand adaptive mesh refinement to maintain solution accuracy.

BEM has emerged as a promising alternative for solving multiple problem. The classical Nystr\"{o}m method~\cite{CK2019Inverse,kress1995numerical,kress1993numerical} can be directly applied to multi-particle configurations. To further enhance computational efficiency, the multiple expansion method (MEM,~\cite{ammari2018field,barucq2018numerical,martin2006multiple}) was developed for problems involving multiple disjoint disks. MEM has since been extended to accommodate arbitrarily shaped particles by enclosing each within an artificial disk, as demonstrated in the S-matrix methods~\cite{lai2015fast,lai2014fast,lai2019framework}, T-matrix methods~\cite{cai1999large,ganesh2012convergence}, and Dirichlet-to-Neumann methods~\cite{acosta2010coupling,grote2004dirichlet}. While these methods perform well for smooth geometries, they face limitations when handling high-curvature features. Furthermore, the solution of large-scale ill-conditioned systems near resonant frequencies results substantial computational costs.

The shape derivative is a powerful tool for analyzing the impact of shape perturbations on scattering phenomena. It has been widely applied in wave scattering contexts such as shape optimization~\cite{bao2014optimal,li2023stochastic,tissot2020optimal}, inverse scattering~\cite{chen2024solving,hettlich1995frechet,hohage1998convergence,roy1997scattering}, and uncertainty quantification~\cite{hao2018computation,harbrecht2013first}. Traditionally, shape derivatives are formulated using the velocity method~\cite{delfour2011shapes,haslinger2003introduction,sokolowski1992introduction} and are primarily developed for simply connected domains. Extending this framework to configurations involving multiple disjoint particles introduces analytical and computational challenges. However, since the multiple scattering problem is formulated via boundary integral equations, the adjoint system—obtained through shape differentiation—naturally retains a boundary integral structure that remains valid for multiple scatterers. A key computational bottleneck lies in efficiently solving this adjoint system. To the best of our knowledge, there seems no prior work in the literature directly addresses the efficient numerical solution of adjoint equations in the context of multiple scattering.

Additionally, the optimal shape design problem constitutes a non-convex optimization challenge, characterized by a complicated objective landscape with multiple local minima, saddle points, and flat regions~\cite{delfour2011shapes,haslinger2003introduction}. The success of the optimization process is highly sensitive to the quality of the initial guess, which must account for (i) the number of nanoparticles, (ii) their geometric configurations, and (iii) their spatial arrangements.  Recently, data-driven approaches—particularly deep learning techniques—have shown significant promise in predicting initial configurations from large datasets~\cite{gao2022artificial,gao2023bayesian,gao2021machine,zhou2023neural}. A major limitation in this field is the absence of comprehensive datasets that systematically correlate multiple nanoparticle geometries with their absorptance.

\subsection{Contributions}

To overcome the aforementioned limitations, we propose a computational framework that integrates the reduced basis method (RBM) with shape-adaptive basis construction and physics-guided initialization. This framework effectively addresses the non-convex optimization challenge of designing broadband absorbers composed of multiple plasmonic nanoparticles, mitigating both computational costs and sensitivity to initial conditions.

We first formulate the forward scattering problem as a boundary integral equation. To circumvent the complexity associated with shape derivatives for multiple particles, we employ a parametric representation of nanoparticle geometries using a finite set of shape parameters. This transforms the original infinite-dimensional shape optimization problem into a finite-dimensional parameter optimization problem constrained by boundary integral equations. The derivative of the objective function with respect to these shape parameters is computed using an adjoint-based approach in conjunction with the shape derivatives of boundary integral operators. For practical implementation, we focus on elliptical nanoparticle geometries. Ellipses offer a favorable balance between geometric complexity and manufacturability: they provide sufficient tunability to modulate plasmonic resonances while remaining simpler to fabricate than more intricate shapes, and more versatile than circular disks.

Next, the high computational cost of broadband simulations necessitates efficient solvers for both the forward and adjoint multiple-scattering problems, as well as careful treatment of singularities in the derivatives of boundary integral operators. Building on the mathematical theory of plasmonic resonances~\cite{ammari2016surface,AFKRYZ2018Mathematical,ammari2017mathematical,ammari2018heat,ammari2016mathematical,yu2018plasmonic}, we leverage the fact that the eigenfunctions of the NP operator form an orthogonal basis in $L^2(\partial D)$ under a suitably defined inner product. Following the MEM in~\cite{fitzpatrick2021convergence}, our RBM expands the boundary density functions of the forward scattering problem in terms of these eigenfunctions, resulting in a semi-discrete formulation. Instead of employing a conventional Galerkin projection, we apply a collocation scheme to obtain a fully discrete system. For the adjoint problem, we utilize the eigenfunctions of the adjoint parameterized NP operator, which substantially improves computational efficiency. Notably, for elliptical particles, both the NP operator and its adjoint admit closed-form eigenfunction representations, enabling exact and efficient computations.

Finally, we develop a physics-informed initialization strategy based on weak scattering approximations. Recognizing that a single elliptical particle supports two resonance modes whose absorptance depend on both the aspect ratio and orientation angle, we construct a off-line dataset of ellipses. Following multiple scattering theory~(\cite{martin2006multiple}), we approximate the total absorptance of weakly interacting particles as the superposition of their individual responses. This formulation leads to a quadratic integer programming problem for generating initial design guesses. Although such problems are generally NP-hard, we adopt a two-stage strategy: first, we relax the integer constraints and solve a continuous quadratic programming problem; next, we obtain an integer solution by rounding the continuous solution and apply a heuristic algorithm to refine it, using the rounded result as an initial search point. This approach provides robust initial estimates for the number of particles, their geometric parameters, and spatial arrangement.

The key features of our optimal design framework include:
\begin{itemize}[leftmargin=*]
	\item A generalized parameterized design formulation for multiple scatterers of arbitrary geometry, naturally extensible to diverse boundary conditions and material systems.
	
	\item A RBM with automatic basis updates during shape optimization, employing distinct tailored bases for the forward and adjoint problems to enhance computational accuracy, especially in the resonance regime.
	
	\item Rigorous mathematical treatment of singularities in boundary integral operators and their shape derivatives by singular splitting, motivated by spectral analysis of the NP operator.
	
	\item A physics-informed initialization strategy based on weak scattering approximations and resonance characteristics specific to plasmonic nanostructures.
	
	\item Demonstrated computational effectiveness through comprehensive numerical experiments, achieving high accuracy even in challenging scenarios involving high-curvature features.
	
\end{itemize}

\subsection{Outline}

The remainder of this paper is structured as follows. \Cref{sec:problem_description} establishes the mathematical formulation of the optimal design problem and introduces key notation. \Cref{sec:parameterized_problem} presents our parameterized reformulation, which transforms the original shape optimization into a finite-dimensional problem. In \Cref{sec:gradient}, we rigorously derive the gradient of the objective function for multi-particle systems. \Cref{sec:RBM} introduces our RBM framework for the efficient solution of both forward and adjoint problems. The complete optimal design algorithm is detailed in \Cref{sec:algorithm}. \Cref{sec:numerical} presents comprehensive numerical results validating the efficiency and accuracy of our approach. Finally, we conclude with a discussion and future research directions in \Cref{sec:conclusion}.

\subsection{Notations and  Preliminary}

Throughout this paper, we denote the wavelength by $\lambda$ and the angular frequency by $\omega$. Let $M$ denote the number of nanoparticles, $P$ the number of parameters for each nanoparticle, and $N$ the number of eigenfunctions used for each particle. We define the parameter space as $W := \mathbb{R}^{MP}$, and let $V := L^2([0, 2\pi), \mathbb{C}^M)$ denote the space of square-integrable, vector-valued functions. The notation $\langle \cdot, \cdot \rangle_X$ indicates the standard inner product in the Hilbert space $X$, while $\|\cdot\|_Y$ represents the associated norm in the normed space $Y$. Subscripts are used to explicitly specify the corresponding function space for each operator.

We define $G^{k}(z) = -\frac{\mathrm{i}}{4} H_0^{(1)}(k|z|)$ as the Green's function for the Helmholtz equation, where $\mathrm{i}$ is the imaginary unit and $H_0^{(1)}$ denotes the zeroth-order Hankel function of the first kind. Next, we define the single layer potential:
\begin{align*}
		\mathcal{S}_D^k[\psi](x) := \int_{\partial D} G^k(x-y) \psi(y)\mathrm{d}\sigma(y), \quad x \in \mathbb{R}^2, 
\end{align*}
where $\psi \in L^2(\partial D)$ is the density function and $\mathrm{d}\sigma$ denotes the surface measure on $\partial D$. The single layer potential exhibits jump relations across the boundary (\cite{AFKRYZ2018Mathematical}):
\begin{align}\label{eq:jump_formula}
	\left. \frac{\partial}{\partial \nu} \mathcal{S}_D^{k}[\psi] \right|_{\pm}(x) = \left( \pm \frac{1}{2} \mathcal{I} + (\mathcal{K}_D^{k})^* \right) [\psi](x), \quad x \in \partial D, 
\end{align}
where $\psi \in L^2(\partial D)$ and the $\pm$ subscripts denote the exterior and interior limits, respectively. The associated boundary integral operator $(\mathcal{K}_D^{k})^*$ is given by:
\begin{align*}
	(\mathcal{K}_D^{k})^*[\psi](x) := \int_{\partial D} \frac{\partial G^{k}(x-y)}{\partial \nu_x} \psi(y) \, \mathrm{d}\sigma(y), \quad x \in \partial D.
\end{align*}
where $\nu_y$ denotes the outward unit normal vector at $y \in \partial D$. Similarly, let $G(z) = \frac{1}{2\pi} \ln |z|$ be the Green's function for the Laplace equation. The single layer potential $\mathcal{S}_D$ and boundary operator $\mathcal{K}_D^*$ for the Laplace equation can be similarly defined with kernel $G(z)$ and also satisfy the jump formula \eqref{eq:jump_formula}.

\section{Problem Description}\label{sec:problem_description}
In what follows, we present a comprehensive mathematical formulation of the multiple scattering and optimal design problems.

\subsection{Boundary integral equation}

We employ layer potential techniques to solve the multiple scattering problem \eqref{eq:original_scattering_problem}, whose solution admits the following integral representation:
\begin{equation}\label{eq:integral_representation}
	u(x) =
	\begin{cases}
		\mathcal{S}_D^{k_c}[\phi](x), & x \in D, \\
		u^i(x) + \mathcal{S}_D^{k_m}[\varphi](x), & x \in \mathbb{R}^2 \setminus \overline{D}.
	\end{cases}
\end{equation}
To satisfy the transmission conditions on $\partial D$ and applying the jump relations \eqref{eq:jump_formula}, we obtain the following boundary integral equation:
\begin{equation}\label{eq:integral_scattering_problem}
	\begin{cases}
		\mathcal{S}_D^{k_c}[\phi] - \mathcal{S}_D^{k_m}[\varphi] = u^i, \\
		\dfrac{1}{\varepsilon_c}\left(-\dfrac{1}{2} \mathcal{I} + (\mathcal{K}_D^{k_c})^*\right)[\phi]
		- \dfrac{1}{\varepsilon_m}\left(\dfrac{1}{2} \mathcal{I} + (\mathcal{K}_D^{k_m})^*\right)[\varphi]
		= \dfrac{1}{\varepsilon_m} \dfrac{\partial u^i}{\partial \nu}.
	\end{cases}
\end{equation}
According to \cite{AFKRYZ2018Mathematical}, this system has a unique solution $(\phi, \varphi) \in L^2(\partial D) \times L^2(\partial D)$.

\subsection{Optimal shape design problem}

We begin by introducing the key concepts of energy flux and energy flow for scalar waves. The Poynting vector is commonly employed to quantify energy flux, and its simplified form for the TM case is presented in the following definition.
\begin{definition}[Energy Flux \cite{born2013principles_of_optics}]
	The energy flux for a scalar wave is given by
	\begin{align*}
		F(x) = -\mathrm{i} C \left[ \overline{u}(x) \nabla u(x) - u(x) \nabla \overline{u}(x) \right],
	\end{align*}
	where $C =  1/ 2\omega\varepsilon_m$ for the transverse magnetic (TM) case. By decomposing the total field as $u = u^i + u^s$, the energy flux can be written as a sum of three components $F = F^i + F^s + F'$, where the incident flux $F^i$, scattered flux $F^s$, and interference flux $F'$ are defined by
	\begin{equation}\label{eq:def_of_Flux}
		\begin{aligned}
			F^i(x) &= -\mathrm{i} C \left[ \overline{u^i}(x) \nabla u^i(x) - u^i(x) \nabla \overline{u^i}(x) \right], \\
			F^s(x) &= -\mathrm{i} C \left[ \overline{u^s}(x) \nabla u^s(x) - u^s(x) \nabla \overline{u^s}(x) \right], \\
			F'(x) &= -\mathrm{i} C \left[ \overline{u^i}(x) \nabla u^s(x) - u^s(x) \nabla \overline{u^i}(x) 
			- u^i(x) \nabla \overline{u^s}(x) + \overline{u^s}(x) \nabla u^i(x) \right].
		\end{aligned}
	\end{equation}
\end{definition}

Next, we consider the energy flow through a partial spherical surface.
\begin{definition}[Energy Flow \cite{born2013principles_of_optics}]
	Let $\partial B_{R,\Theta}$ denote the partial sphere
	\begin{align}\label{eq:partial_sphere}
		\partial B_{R,\Theta} = \left\{ x \in \mathbb{R}^2\mid x = R(\cos\theta, \sin\theta),~\theta\in \Theta \right\},~\text{where}~\Theta=(\bar{\theta} - \Delta\theta, \bar{\theta} + \Delta\theta),
	\end{align}
	illustrated in \Cref{fig:multiple_nanoparticles}. The total energy flow through this surface is defined by
	\begin{align}\label{eq:def_of_Flow}
		E_{R,\Theta} = \int_{\partial B_{R,\Theta}} F(x) \cdot \nu(x)\, \mathrm{d}\sigma(x),
	\end{align}
	where $F$ is the energy flux and $\nu$ is the outward unit normal. We decompose the energy into incident, scattered, and interference parts, i.e., $E_{R,\Theta} = E^i_{R,\Theta} + E^s_{R,\Theta} + E'_{R,\Theta}$, where each component corresponds to the integral of $F^i$, $F^s$, and $F'$, respectively. Consequently, the absorbed energy satisfies $E^{a}_{R,\Theta} = E^i_{R,\Theta} - E_{R,\Theta} = -E^s_{R,\Theta} - E'_{R,\Theta}$.
\end{definition}

To facilitate the subsequent optimal design problem, we next derive an asymptotic representation for these integrals. The following result, proved in \cref{appendix:proof_of_energy_flow}, provides the leading-order behavior as $R \to \infty$.
\begin{theorem}[Asymptotic of Energy Flow]\label{thm:asymptotic_of_energy_flow}
	The energy flow have the following asymptotic behavior as $R \to \infty$
	\begin{align*}
			E^i_{R,\Theta} = 4C k_m R \sin(\Delta\theta) \cos(\bar{\theta} - \theta_0), \
			E^s_{R,\Theta} = 2C k_m  \| u^\infty(\hat{x}(\theta)) \|^2_{L^2(\Theta)}  + \mathcal{O}(R^{-1}),
	\end{align*}
	where $\hat{x}(\theta) = (\cos\theta, \sin\theta)$. The asymptotic behavior of $E^\prime_{R,\Theta}$ exhibits two cases:
	\begin{itemize}[leftmargin=*]
		\item Forward Scattering ($\theta_0 \in \Theta$): $E^\prime_{R,\Theta} = 2C k_m \sqrt{\frac{8\pi}{k_m}} \Im\left(e^{\mathrm{i} \frac{3\pi}{4}} u^\infty(d)\right) + \mathcal{O}\left(R^{-1/2}\right)$.
		\item Backward Scattering ($\theta_0 + \pi \in \Theta$): $E^\prime_{R,\Theta} = \mathcal{O}\left(R^{-1/2}\right)$.
	\end{itemize}
	Here $\theta_{0}$ is the incident angle and the far‐field pattern $u^{\infty}$ is defined by
	\begin{equation}\label{eq:far_field_operator}
		u^\infty(\hat{x}) = -\frac{e^{\mathrm{i} \pi /4}}{\sqrt{8\pi k_m}} \int_{\partial D} e^{-\mathrm{i} k_m \hat{x} \cdot y} \varphi(y) \, \mathrm{d}\sigma(y).
	\end{equation}
\end{theorem}

\begin{proposition}[Optical Theorem]\label{prop:optical_theorem}
	For $\Delta\theta = \pi$, we recover the two-dimensional optical theorem:
	\begin{align}\label{eq:cross_section}
		Q^e = -\frac{E^\prime}{|F^i|} = \sqrt{\frac{8\pi}{k_m}} \Im\left[e^{\mathrm{i} \frac{3\pi}{4}} u^\infty(d)\right],\quad
		Q^s = \frac{E^s}{|F^i|} = \int_0^{2\pi} \left| u^\infty(\hat{x}(\theta)) \right|^2 \, \mathrm{d}\theta,
	\end{align}
	where $Q^{e}$ is the extinction cross-section, expressed in terms of the imaginary part of the forward-scattering amplitude multiplies a phase, $Q^s$ is the scattering cross-section, given by the far-field intensity integrated over the unit circle $\mathbb{S}^1$, and $Q^a = Q^e - Q^s$ is the absorption cross-section. Note that the expression for  $Q^{e}$ in $\mathbb{R}^2$ differs fundamentally from its three-dimensional analogue in $\mathbb{R}^3$ (cf.\cite{ammari2017mathematical} for comparison).
\end{proposition}

\begin{remark}
    As established in \cref{thm:asymptotic_of_energy_flow} and the Optical Theorem (\cref{prop:optical_theorem}), total extinction is governed by the imaginary part of the forward-scattering amplitude. This is due to the distinctive feature of forward scattering: only in this direction do the incident and scattered waves interfere persistently, resulting in energy loss along the incident path. The Optical Theorem formalizes this connection by directly relating extinction to forward scattering. Consequently, the forward direction is primarily responsible for changes in the energy flux. For these reasons, our subsequent analysis focuses on forward scattering phenomena.
\end{remark}

Based on the preceding analysis and preparation, we define the \emph{energy absorptance} on the partial sphere $\partial B_{R,\Theta}$ as
\begin{align}\label{eq:absorptance}
	A(\varphi, D, \lambda) :=E^a_{R,\Theta}/E^i_{R,\Theta}
	= -(E^s_{R,\Theta} + E^\prime_{R,\Theta})/E^i_{R,\Theta},
\end{align}
where $A$ depends on both the wavelength $\lambda$ and the nanoparticle configuration $D$. Let $A^{\text{tar}}(\lambda)$ denote the target absorptance spectrum over the wavelength range $\Lambda := [\lambda_{\min}, \lambda_{\max}]$. We formulate the objective functional:
\begin{equation}\label{eq:objective_functional}
	J(\phi, \varphi, D, \lambda) = \|A(\varphi, D, \lambda) - A^{\mathrm{tar}}(\lambda)\|^2_{L^2(\Lambda)}
	= \int_{\Lambda} \left| A(\varphi, D, \lambda) - A^{\mathrm{tar}}(\lambda) \right|^2 \, \mathrm{d}\lambda.
\end{equation}
The optimal design problem is then formulated as the constrained minimization:
\begin{equation}\label{eq:infinite_optimal_design}
	\begin{aligned}
		\min_{(\phi, \varphi, D, \lambda) \in U_{\text{ad}}} \ & J(\phi, \varphi, D, \lambda) 
		\quad \text{subject to} \ \eqref{eq:integral_scattering_problem},
	\end{aligned}
\end{equation}
where $U_{\text{ad}}$ denotes the set of admissible designs.

\section{Parameterization of the Design Problem}\label{sec:parameterized_problem}

In this section, we first consider the parameterized form of multiple particles. Next, we introduce the parameterized form of the boundary operator. Finally, we obtain the parameterized optimal design problem, which is a reduced finite-dimensional form of the original infinite-dimensional optimal design problem \eqref{eq:infinite_optimal_design}. 

\subsection{Parameterized multiple particles}

We assume that each nanoparticle has an analytic boundary, represented by a $2\pi$-periodic parametric curve. Let $w_m \in \mathbb{R}^P$ parametrize the boundary 
\begin{align*}
	\partial D_m := \left\{ x(t; w_m) = (x_1(t; w_m), x_2(t; w_m)) \in \mathbb{R}^2 \,\middle|\, t \in [0, 2\pi) \right\},
\end{align*}
where $t \mapsto x(t; w_m)$ is analytic, $2\pi$-periodic, counterclockwise-oriented, and satisfies $|x'(t; w_m)| > 0$ for all $t$. For points $x(t; w_n) \in \partial D_n$ and $y(s; w_m) \in \partial D_m$, the relative displacement is $z(t, s; w_n, w_m) = x(t; w_n) - y(s; w_m)$ for $n \neq m$ and $z(t, s; w_m) = x(t; w_m) - y(s; w_m)$ for $n = m$. For a system of $M$ nanoparticles with parameters $w = (w_1, \dots, w_M)$, the total boundary is $\partial D = \{x(t; w) \mid t \in [0, 2\pi)\}$, where $x(t; w) = \big(x(t; w_1), \dots, x(t; w_M)\big)^\top$. The system-wide relative displacement is $z(t, s; w) = x(t; w) - y(s; w)^\top$.

\subsection{Parameterized boundary operators}

Based on the boundary parametrization, we define the following boundary integral operators.

\begin{definition}[Parametrized Boundary Operators]\label{def:parameterized_BDO}
	Let $\tilde{\psi}_m(t) := \psi(y(t; w_m))$ denote the parametrized density function on the particle $D_m$. The corresponding boundary integral operator evaluated on the particle $D_m$ is defined as follows:
	\begin{align*}
	\mathcal{S}_{w_n, w_m}^k[\tilde{\psi}_m](t) := \int_{0}^{2\pi} G^k(z(t, s; w_n, w_m)) \, |y^{\prime}(s; w_m)| \, \tilde{\psi}_m(s) \, \mathrm{d}s.
	\end{align*}
	For the full system with $\tilde{\psi}(t) := (\tilde{\psi}_1(t), \dots, \tilde{\psi}_M(t))^\top\in V$, the parametrized single-layer boundary operator on $\partial D$ is defined as the matrix form
	$\mathcal{S}_{w}^k=[\mathcal{S}_{w_n, w_m}^k]$. Similarly, the parametrized boundary operator $\left(\mathcal{K}_w^k\right)^*=[\left(\mathcal{K}_{w_n, w_m}^k\right)^*]$ where each component operator is defined by
	\begin{align*}
			\left(\mathcal{K}_{w_n, w_m}^k\right)^*[\tilde{\psi}_m](t) &:= \int_{0}^{2\pi} \frac{\partial G^k(z(t, s; w_n, w_m))}{\partial \nu_x(t; w_n)} \, |y^{\prime}(s; w_m)| \, \tilde{\psi}_m(s) \, \mathrm{d}s.
	\end{align*}
	For the diagonal term $m=n$,  denote $\mathcal{S}_{w_m}^k$ and $\left(\mathcal{K}_{w_m}^k\right)^*$ for convenience.
\end{definition}

Next, we give the definition of adjoint operators of the parametrized boundary operators $\mathcal{S}_w^k$ and $(\mathcal{K}_w^k)^*$.
\begin{definition}[Adjoint Parametrized Operators]\label{def:adjoint_operator}
	The adjoint operators of $\mathcal{S}_w^k$ and $(\mathcal{K}_w^k)^*$ under the inner product $L^2([0, 2\pi), \mathbb{C}^{M})$ are given by
	\begin{equation}\label{eq:Sskw}
		\begin{aligned}
			(\mathcal{S}_w^k)^*[\tilde{\psi}](t) &:= \int_{0}^{2\pi} G^{-k}(z(t, s; w)) \, |x^{\prime}(t; w)| \, \tilde{\psi}(s) \, \mathrm{d}s, \\
			\mathcal{K}_w^k[\tilde{\psi}](t) &:= \int_{0}^{2\pi} \frac{\partial G^{-k}(z(t, s; w))}{\partial \nu_y(s; w)} \, |x^{\prime}(t; w)| \, \tilde{\psi}(s) \, \mathrm{d}s.
		\end{aligned}
	\end{equation}
\end{definition}

\begin{remark}
	The parametrized form of the adjoint operator $(\mathcal{S}_D^{k})^{*}$ is different from that of $(\mathcal{S}_w^{k})^{*}$ (the adjoint of the parametrized operator $\mathcal{S}_{w}^k$). Specifically,
	\begin{align*}
		(\mathcal{S}_D^k)^* \psi = \int_{\partial D} G^{-k}(x(t) - y) \psi(y) \mathrm{d}\sigma(y)
		= \int_0^{2\pi} G^{-k}(z(t, s; w)) |y'(s; w)| \tilde{\psi}(s)\mathrm{d}s.
	\end{align*}
	Similarly, the parametrized adjoint operator $(\mathcal{K}_D^{k})^{*}$ is different from $(\mathcal{K}_w^{k})^{*}$.
\end{remark}

\subsection{Parameterized optimal design problem}

Building upon the parameterization of multiple particles and boundary operators, we are now ready to reformulate the boundary integral system \eqref{eq:integral_scattering_problem} into a parametric form:
\begin{equation}
\resizebox{0.9\textwidth}{!}{$
  \begin{aligned}
    \begin{cases}
      \displaystyle \mathcal{E}[\tilde{\phi}, \tilde{\varphi}, w, \lambda] :=
      \mathcal{S}_w^{k_c}[\tilde{\phi}]
      - \mathcal{S}_w^{k_m}[\tilde{\varphi}]
      - f_1 = 0, \\
     \displaystyle \mathcal{F}[\tilde{\phi}, \tilde{\varphi}, w, \lambda] :=
      \frac{1}{\varepsilon_c} \left( -\frac{1}{2} \mathcal{I} + (\mathcal{K}_w^{k_c})^* \right)[\tilde{\phi}]
      - \frac{1}{\varepsilon_m} \left( \frac{1}{2} \mathcal{I} + (\mathcal{K}_w^{k_m})^* \right)[\tilde{\varphi}]
      - f_2 = 0,
    \end{cases}
  \end{aligned}
$}
\label{eq:parameterized_state_equation}
\end{equation}
where  $f_1 = f_1[w,\lambda] = u^i(x(t; w))$ and $f_2 = f_2[w,\lambda] = \frac{1}{\varepsilon_m} \frac{\partial u^i(x(t; w))}{\partial \nu_x(t; w)}$. These operators define  mappings as: $\mathcal{E}, \mathcal{F}: V \times V \times W \times \mathbb{R} \rightarrow V$ and $f_1, f_2: W \times \mathbb{R} \rightarrow V$.

According to the far-field representation formula \eqref{eq:far_field_operator}, the far-field pattern can be rewritten as
\begin{align*}
	u^{\infty}[\tilde{\varphi}, w, \lambda](\theta) = \frac{-e^{\mathrm{i} \frac{\pi}{4}}}{\sqrt{8 k_m \pi}} \int_{0}^{2\pi} e^{-\mathrm{i} k_m \hat{x}(\theta) \cdot y(s; w)} \, |y^{\prime}(s; w)| \, \tilde{\varphi}(s) \, \mathrm{d}s = \frac{-e^{\mathrm{i} \frac{\pi}{4}}}{\sqrt{8 k_m \pi}} h[\tilde{\varphi}, w, \lambda](\theta),
\end{align*}
where $h[\tilde{\varphi}, w, \lambda](\theta) = \langle \tilde{\varphi}(s), p[w, \lambda](\theta, s) \rangle_{V}$ and  $p[w, \lambda](\theta, s) = e^{\mathrm{i} k_m \hat{x}(\theta) \cdot y(s; w)} \, |y^{\prime}(s; w)|$.

Restricting our analysis to the leading-order term in \cref{thm:asymptotic_of_energy_flow}, we obtain the following expressions for the incident and scattered energy
\begin{equation*}
	\begin{aligned}
		E^i_{R,\Theta}(\lambda) = 4C k_m R \sin(\Delta\theta) \cos(\overline{\theta} - \theta_0),~ E^s_{R,\Theta}(\tilde{\varphi}, w, \lambda) &= \frac{2C k_m}{8 k_m \pi} \| h[\tilde{\varphi}, w, \lambda] \|^2_{L^2(\Theta)}.
	\end{aligned}
\end{equation*}
Moreover, $E^\prime_{R,\Theta}(\tilde{\varphi}, w, \lambda) = 0$ if $\theta_0 \notin \Theta$; otherwise,
\begin{equation*}
	\begin{aligned}
		E^\prime_{R,\Theta}(\tilde{\varphi}, w, \lambda)
		= 2C k_m \Im\left( \frac{1}{k_m} h[\tilde{\varphi}, w, \lambda](\theta_0) \right) = 2C k_m \Re \left( \frac{-\mathrm{i}}{k_m} h[\tilde{\varphi}, w, \lambda](\theta_0) \right).
	\end{aligned}
\end{equation*}
Thus, the absorptance function $A$ is 
\begin{equation*}
	\begin{aligned}
		A(\tilde{\varphi}, w, \lambda) = &-\left(E^s_{R,\Theta}(\tilde{\varphi}, w, \lambda) + E^\prime_{R,\Theta}(\tilde{\varphi}, w, \lambda)\right)/E^i_{R,\Theta}(\lambda)\\
		= & -\frac{1}{L_{R,\Theta}}\left(  \frac{1}{8 k_m \pi} \| h[\tilde{\varphi}, w, \lambda] \|^2_{L^2(\Theta)}+\Re\left( \frac{-\mathrm{i}}{k_m} h[\tilde{\varphi}, w, \lambda](\theta_0) \right)\right),
	\end{aligned}
\end{equation*}
where the normalization constant $L_{R,\Theta} = 2R \sin(\Delta\theta) \cos(\overline{\theta} - \theta_0)$. Thus, the parameterized objective functional is defined by
\begin{equation}\label{eq:parameterized_objective_function}
	J[\tilde{\phi}, \tilde{\varphi}, w, \lambda] := \| A(\tilde{\varphi}, w, \lambda) - A^{\mathrm{tar}}(\lambda) \|^2_{L^2(\Lambda)}.
\end{equation}
We now formulate the optimal design problem in a parameterized form:
\begin{equation}\label{eq:parameterized_optimal_design}
	\begin{aligned}
		\min_{(\tilde{\phi}, \tilde{\varphi}, w, \lambda) \in U_{\mathrm{ad}}} \ & J(\tilde{\phi}, \tilde{\varphi}, w, \lambda)
		\quad \text{subject to} \ &
		\begin{cases}
			\mathcal{E}[\tilde{\phi}, \tilde{\varphi}, w, \lambda] = 0, \\
			\mathcal{F}[\tilde{\phi}, \tilde{\varphi}, w, \lambda] = 0,
		\end{cases}
	\end{aligned}
\end{equation}
where the admissible set is $U_{\mathrm{ad}} = V \times V \times W_{\mathrm{ad}} \times \Lambda$ and $W_{\mathrm{ad}} \subset W$ is a bounded set.

\section{Gradient of Objective Function}\label{sec:gradient}

In this section, we aim to obtain the gradient of the objective function based on the parameterized optimal design problem \eqref{eq:parameterized_optimal_design}, thereby avoiding the calculation of complex shape derivatives for multiple particles. The gradient of the objective functional is derived using the adjoint method.

\subsection{Adjoint method for gradient}

For every  $w$ and $\lambda$, there exists a unique solution $(\tilde{\phi}, \tilde{\varphi})$ to the scattering problem~\eqref{eq:parameterized_state_equation}, where $\tilde{\phi} = \tilde{\phi}(w, \lambda)$ and $\tilde{\varphi} = \tilde{\varphi}(w, \lambda)$. To facilitate the discussion of gradient computation, we introduce the reduced objective function corresponding to~\eqref{eq:parameterized_objective_function}:
\begin{equation*}
	\mathcal{J}(w) := \left\| A(\tilde{\varphi}(w, \lambda), w, \lambda) - A^{\mathrm{tar}}(\lambda) \right\|^2_{L^2(\Lambda)}.
\end{equation*}
The computation of the gradient of $\mathcal{J}(w)$ is established in the following theorem, whose proof is provided in \cref{appendix:gradient}.
\begin{theorem}\label{thm:gradient}
	The gradient of the objective functional $\mathcal{J}(w)$ admits the representation
	\begin{equation}\label{eq:full_gradient}
    \resizebox{0.91\textwidth}{!}{$
		\mathcal{J}^\prime(w) = 2 \left\langle A(\tilde{\varphi}, w, \lambda) - A^{\mathrm{tar}}(\lambda), \ A_w(\tilde{\varphi}, w, \lambda) - \Re \left( \langle \tilde{p}, \mathcal{E}_w[\tilde{\phi}, \tilde{\varphi}, w, \lambda] \rangle_{V} + \langle \tilde{q}, \mathcal{F}_w[\tilde{\phi}, \tilde{\varphi}, w, \lambda] \rangle_{V} \right) \right\rangle_{L^2(\Lambda)},
        $}
	\end{equation}
	where the adjoint densities $\tilde{p}$ and $\tilde{q}$ are determined by the following system:
	\begin{equation}\label{eq:abstract_adjoint_equation}
		\begin{cases}
			\mathcal{E}^*_{\tilde{\phi}}[\tilde{\phi}, \tilde{\varphi}, w, \lambda] \tilde{p} + \mathcal{F}^*_{\tilde{\phi}}[\tilde{\phi}, \tilde{\varphi}, w, \lambda] \tilde{q} = A_{\tilde{\phi}}(\tilde{\varphi}, w, \lambda), \\
			\mathcal{E}^*_{\tilde{\varphi}}[\tilde{\phi}, \tilde{\varphi}, w, \lambda] \tilde{p} + \mathcal{F}^*_{\tilde{\varphi}}[\tilde{\phi}, \tilde{\varphi}, w, \lambda] \tilde{q} = A_{\tilde{\varphi}}(\tilde{\varphi}, w, \lambda).
		\end{cases}
	\end{equation}
\end{theorem}

\subsection{Explicit form of gradient terms}

In this subsection, we further derive the explicit expressions for the terms appearing in the full gradient $\mathcal{J}^\prime(w)$ of the objective functional as described in \eqref{eq:full_gradient}. These terms include the adjoint densities $\tilde{p}$ and $\tilde{q}$, as well as the gradients of operators $A_w(\tilde{\varphi}, w, \lambda)$, $\mathcal{E}_w[\tilde{\phi}, \tilde{\varphi}, w, \lambda]$, and $\mathcal{F}_w[\tilde{\phi}, \tilde{\varphi}, w, \lambda]$.

We begin by considering the abstract adjoint equation \eqref{eq:abstract_adjoint_equation}. The adjoint operators on the left-hand side are explicitly given by:
\begin{align*}
	\mathcal{E}^*_{\tilde{\phi}}(\tilde{\phi},\tilde{\varphi}, w)
		&= (\mathcal{S}_w^{k_c})^*,
		& \mathcal{F}^*_{\tilde{\phi}}(\tilde{\phi},\tilde{\varphi}, w)
		&= (-\frac{1}{2} \mathcal{I}+\mathcal{K}_w^{k_c})/\overline{\varepsilon_c}, \\
		\mathcal{E}^*_{\tilde{\varphi}}(\tilde{\phi},\tilde{\varphi}, w)
		&= -(\mathcal{S}_w^{k_m})^*,
		& \mathcal{F}^*_{\tilde{\varphi}}(\tilde{\phi},\tilde{\varphi}, w)
		&= -(\frac{1}{2} \mathcal{I}+\mathcal{K}_w^{k_m})/\varepsilon_m,
\end{align*}
where the adjoint operators $(\mathcal{S}_w^{k})^*$ and $\mathcal{K}_w^{k}$ have been defined in \eqref{eq:Sskw}. Substituting these into \eqref{eq:abstract_adjoint_equation}, the adjoint system becomes:
\begin{equation}\label{eq:adjoint_equation}
\begin{aligned}
    	\begin{cases}
		 \displaystyle (\mathcal{S}_w^{k_c})^*[\tilde{p}] + \frac{1}{\overline{\varepsilon_c}} \left( -\frac{1}{2} \mathcal{I} + \mathcal{K}_w^{k_c} \right)[\tilde{q}] = g_1, \\
		\displaystyle -(\mathcal{S}_w^{k_m})^*[\tilde{p}] - \frac{1}{\varepsilon_m} \left( \frac{1}{2} \mathcal{I} + \mathcal{K}_w^{k_m} \right)[\tilde{q}] = g_2,
	\end{cases}
\end{aligned}
\end{equation}
where $g_1 = A_{\tilde{\phi}}(\tilde{\phi}, \tilde{\varphi}, w, \lambda)$ and $g_2 = A_{\tilde{\varphi}}(\tilde{\varphi}, w, \lambda)$. Notably, $A_{\tilde{\phi}}(\tilde{\varphi}, w, \lambda) = 0$ due to the independence of $A(\tilde{\varphi}, w, \lambda)$ from $\tilde{\phi}$. Therefore, we focus on the explicit form of $A_{\tilde{\varphi}}(\tilde{\varphi}, w, \lambda)$, which can be decomposed as:
\begin{align*}
	A_{\tilde{\varphi}}(\tilde{\varphi}, w, \lambda) =& -\frac{1}{L_{R,\Theta}}\left( \frac{1}{4 k_m \pi} \left\langle h[\tilde{\varphi}, w, \lambda], \overline{p[w, \lambda]} \right\rangle_{L^2(\Theta)}+ \frac{\mathrm{i}}{k_m} p[w, \lambda]\right).
\end{align*}

Then, we calculate each term that arises in the expression for gradient of operators $A_w(\tilde{\varphi}, w, \lambda)$,  $\mathcal{E}_w[\tilde{\phi}, \tilde{\varphi}, w, \lambda]$, and $ \mathcal{F}_w[\tilde{\phi}, \tilde{\varphi}, w, \lambda]$. First, we consider
\begin{equation}\label{eq:dAdw}
	\begin{aligned}
		A_{w}(\tilde{\varphi}, w, \lambda) = \frac{-1}{L_{R, \Theta}} \Re \left(\frac{1}{4k_m\pi}  \left\langle h_w[\tilde{\varphi}, w, \lambda], h[\tilde{\varphi}, w, \lambda] \right\rangle_{L^2(\Theta)} +  \frac{\mathrm{i}}{k_m} h_w[\tilde{\varphi}, w, \lambda] \right), 
	\end{aligned}
\end{equation}
where $h_w[\tilde{\varphi}, w, \lambda](\theta) = \left\langle \tilde{\varphi}(s), \, p_w[w, \lambda](\theta, s) \right\rangle_{V}$ and $p_w[w, \lambda](\theta, s)$ satisfies:
\begin{align*}
	p_w[w, \lambda](\theta, s) = \left(\mathrm{i} k_m \, \hat{x}(\theta) \cdot \frac{d y(s; w)}{dw} + \frac{1}{|y'(s; w)|} \frac{d |y'(s; w)|}{dw} \right) p[w, \lambda](\theta, s).
\end{align*}

Second, we consider the vector-valued functions $\mathcal{E}_w(\tilde{\phi}, \tilde{\varphi}, w, \lambda)$ and $\mathcal{F}_w(\tilde{\phi}, \tilde{\varphi}, w, \lambda)$, which describe the derivatives of the operators with respect to the shape parameter $w$. They are given by:
\begin{equation}\label{eq:dEdw_and_dFdw}
\begin{aligned}
    	\begin{cases}
		\displaystyle	\mathcal{E}_w(\tilde{\phi}, \tilde{\varphi}, w, \lambda) = \frac{\partial \mathcal{S}_w^{k_c}}{\partial w}[\tilde{\phi}] - \frac{\partial \mathcal{S}_w^{k_m}}{\partial w}[\tilde{\varphi}] - \frac{\partial f_1[w, \lambda]}{\partial w}, \\
		\displaystyle	\mathcal{F}_w(\tilde{\phi}, \tilde{\varphi}, w, \lambda) = \frac{1}{\varepsilon_c} \frac{\partial (\mathcal{K}_w^{k_c})^*}{\partial w}[\tilde{\phi}] - \frac{1}{\varepsilon_m} \frac{\partial (\mathcal{K}_w^{k_m})^*}{\partial w}[\tilde{\varphi}] - \frac{\partial f_2[w, \lambda]}{\partial w}.
	\end{cases}
\end{aligned}
\end{equation}
In order to compute the derivatives of $f_1[w, \lambda]$ and $f_2[w, \lambda]$, we use the fact that the incident wave is given by $u^i(x) = e^{\mathrm{i}k_m d \cdot x}$. Then, we have
\begin{align*}
	f_1[w, \lambda] = &  \left(u^i(x(t; w_1)), \cdots, u^i(x(t; w_M))\right)^{\top},\\
	f_2[w, \lambda] = &  \left(\mathrm{i}k_m d \cdot \nu(t; w_1) u^i(x(t; w_1)), \cdots, \mathrm{i}k_m d \cdot \nu(t; w_M) u^i(x(t; w_M))\right)^{\top}/\varepsilon_m.
\end{align*}
We denote the derivatives with respect to $w$ component-wise as:
\begin{equation*}
	\begin{aligned}
		\frac{\partial f_1}{\partial w} =& (\frac{\partial f_1}{\partial w_1}, \cdots, \frac{\partial f_1}{\partial w_M})^\top,&
		\frac{\partial f_2}{\partial w} =& (\frac{\partial f_2}{\partial w_1}, \cdots, \frac{\partial f_2}{\partial w_M})^\top \in L^2([0, 2\pi),W).
	\end{aligned}
\end{equation*}
Then, for each $n = 1, \dots, M$, $\frac{\partial f_1[w,\lambda]}{\partial w_n}$ and $\frac{\partial f_2[w,\lambda]}{\partial w_n}$ can be computed explicitly.

Before proceeding further, we provide the definition of the derivative of the boundary operators with respect to the shape parameter $w$.
\begin{definition}
	The derivatives of the parameterized boundary layer potentials with respect to shape parameters are given by
	\begin{equation}\label{eq:derivative_of_boundary_op}
		\begin{aligned}
			\frac{\partial \mathcal{S}_{w_n,w_m}^{k}}{\partial w_p}[\tilde{\psi}_m] &:= \int_{0}^{2\pi} \left[ \frac{\partial}{\partial w_p} \left( G^k(z(t,s;w_n,w_m)) \, |y'(s;w_m)| \right) \right] \tilde{\psi}_m(s) \, \mathrm{d}s, \\
			\frac{\partial (\mathcal{K}_{w_n,w_m}^{k})^*}{\partial w_p}[\tilde{\psi}_m] &:= \int_{0}^{2\pi} \left[ \frac{\partial}{\partial w_p} \left( \frac{\partial G^k(z(t,s;w_n,w_m))}{\partial \nu_x(t;w_n)} \, |y'(s;w_m)| \right) \right] \tilde{\psi}_m(s) \, \mathrm{d}s,
		\end{aligned}
	\end{equation}
	where $1 \leq m,n,p \leq M$. Similarly, if the shape parameter $w_p$ affects the density $\tilde{\psi}_m$ as well, then the full derivative of the boundary operator applied to $\tilde{\psi}_m$ is defined as $\frac{\partial \mathcal{S}_{w_n,w_m}^{k}[\tilde{\psi}_m]}{\partial w_p}$ and $\frac{\partial (\mathcal{K}_{w_n,w_m}^{k})^*[\tilde{\psi}_m]}{\partial w_p}$. These derivatives are non-zero only when $p = m$ or $p = n$. Moreover, if the density function $\tilde{\psi}_m$ is independent of $w_p$, then the two notions of derivatives coincide: $\frac{\partial \mathcal{S}_{w_n,w_m}^{k}}{\partial w_p}[\tilde{\psi}_m] = \frac{\partial \mathcal{S}_{w_n,w_m}^{k}[\tilde{\psi}_m]}{\partial w_p}$ and $\frac{\partial (\mathcal{K}_{w_n,w_m}^{k})^*}{\partial w_p}[\tilde{\psi}_m] = \frac{\partial (\mathcal{K}_{w_n,w_m}^{k})^*[\tilde{\psi}_m]}{\partial w_p}$. The derivative of $\mathcal{S}_w^k$ with respect to $w$ is defined by
	\begin{equation*}
		\frac{\partial \mathcal{S}_w^k}{\partial w}[\tilde{\psi}] := \left( \frac{\partial \mathcal{S}_w^k}{\partial w_1}[\tilde{\psi}], \cdots, \frac{\partial \mathcal{S}_w^k}{\partial w_M}[\tilde{\psi}] \right)^\top \in L^2([0, 2\pi),W),
	\end{equation*}
	where the $m$-th component is given by
	\begin{equation}\label{eq:dSkwdwn}
		\begin{aligned}
			\frac{\partial \mathcal{S}_w^k}{\partial w_m}[\tilde{\psi}]
			=
			\begin{bmatrix}
				0 & \cdots & \frac{\partial \mathcal{S}_{w_1,w_m}^k}{\partial w_m} & \cdots & 0 \\
				\vdots & \cdots & \vdots & \cdots & \vdots \\
				\frac{\partial \mathcal{S}_{w_m,w_1}^k}{\partial w_m} & \cdots & \frac{\partial \mathcal{S}_{w_m}^k}{\partial w_m} & \cdots & \frac{\partial \mathcal{S}_{w_m,w_M}^k}{\partial w_m} \\
				\vdots & \cdots & \vdots & \cdots & \vdots \\
				0 & \cdots & \frac{\partial \mathcal{S}_{w_M,w_m}^k}{\partial w_m} & \cdots & 0
			\end{bmatrix}
			\begin{bmatrix}
				\tilde{\psi}_1 \\
				\vdots \\
				\tilde{\psi}_m \\
				\vdots \\
				\tilde{\psi}_M
			\end{bmatrix}.
		\end{aligned}
	\end{equation}
	The derivative of $(\mathcal{K}_w^k)^*$ are analogous and therefore omitted for brevity.
\end{definition}

\section{Reduced Basis Method}\label{sec:RBM}

In this section, we introduce the RBM for solving the scattering problem and its adjoint, as well as for computing derivatives of boundary integral operators. The RBM uses adaptive, shape-dependent basis functions, effectively mitigating the singular behavior of boundary layer potentials. It also maintains high accuracy in resonant regimes and for high-curvature geometries.

\subsection{Spectral of the NP operator and its adjoint}

We begin by reviewing the spectral properties of the Neumann–Poincaré (NP) operator $\mathcal{K}_D^*$, which characterizes plasmon resonances in the quasi-static regime. The spectral properties of $\mathcal{K}_D^*$ have recently attracted significant attention due to their applications in plasmonics \cite{ammari2016surface, ammari2017mathematical}. As shown in \cite{AFKRYZ2018Mathematical}, the eigenfunctions of $\mathcal{K}_D^*$ form a complete basis for $L^2(\partial D)$, enabling the expansion of boundary layer densities in terms of these eigenfunctions. For certain canonical geometries—including disks, ellipses, concentric disks, and confocal ellipses—explicit expressions for both the eigenvalues and eigenfunctions of $\mathcal{K}_D^*$ are available. These explicit forms are particularly useful for the theoretical analysis and numerical approximation of plasmonic resonances.

Next, we focus on elliptical particle shapes within the RBM framework. This geometric choice offers several key advantages. First, elliptical domains allow for explicit spectral characterization, as the complete spectrum and eigenfunctions of the NP operator $\mathcal{K}_D^*$ are known analytically (see \cite{AFKRYZ2018Mathematical}), enabling efficient and accurate computation. Second, ellipses strike a favorable balance between geometric flexibility and computational efficiency. They are parameterized by only two independent variables—the semi-axes $a$ and $b$—along with  translations and rotations. This provides greater versatility than circular shapes without significantly increasing computational complexity. Third, ellipses serve as a unified shape representation, naturally interpolating between several canonical geometries. For instance, they reduce to circles when $a = b$ and approximate slender, quasi-rectangular shapes when one axis is much larger than the other (i.e., $a \gg b$ ). Finally, elliptical particles are capable of exhibiting broadband resonance behaviors and are easy to construct in practical applications.

Then, we introduce elliptic coordinates and review the spectral properties of the NP operator $\mathcal{K}_D^*$. For a point $x = (x_1, x_2)$ in Cartesian coordinates, the elliptic coordinates $(\rho, t)$ are defined by  
\begin{equation}\label{eq:elliptic_coordinate}
	x_1 = c \cosh \rho \cos t, \quad
	x_2 = c \sinh \rho \sin t, \quad
	\rho > 0,~ 0 \leq t < 2\pi,
\end{equation}
where $c > 0$ is the focal distance. The level curve defined by $\{x = (x_1, x_2) \mid \rho = \rho_0,~ 0 \leq t < 2\pi\}$ corresponds to an ellipse with foci located at $(\pm c, 0)$. The geometric parameters of the ellipse are related to $\rho_0$ via the identities  $a = c \cosh \rho_0$, $b = c \sinh \rho_0$, and $\rho_0 = \ln(a + b) - \ln c$, where $a$ and $b$ denote the semi-major and semi-minor axes, respectively. The differential geometric quantities on the ellipse are expressed as  
\begin{equation}\label{eq:differential_property_of_ellipse}
	d\sigma = \Xi(\rho, t) \, dt, \quad
	\partial_\nu = 1/\Xi(\rho, t) \partial_\rho, \quad
	\Xi(\rho, t) = c \sqrt{\sinh^2 \rho + \sin^2 t},
\end{equation}
where $\Xi(\rho, t)$ is the metric coefficient associated with the elliptic coordinate system.

\begin{definition}[Elliptic Particle]\label{def:elliptic_particle}
	Let $w_m = (a_m, b_m, \theta_m, x_{1,m}, x_{2,m})$ denote the parameterization of the $m$-th elliptical particle $D_m$, where $a_m$ and $b_m$ are the semi-major and semi-minor axes; $\theta_m$ is the counterclockwise rotation angle; $(x_{1,m}, x_{2,m})$ is the center of the particle. The boundary $\partial D_m$ in Cartesian coordinates is represented by the parametric form:
	\begin{equation}\label{eq:parameterized_ellipse}
		x(t; w_m)=
		\begin{bmatrix}
			x_1(t; w_m) \\
			x_2(t; w_m)
		\end{bmatrix}
		=
		\begin{bmatrix}
			\cos \theta_m & -\sin \theta_m \\
			\sin \theta_m & \cos \theta_m
		\end{bmatrix}
		\begin{bmatrix}
			a_m \cos t \\
			b_m \sin t
		\end{bmatrix}
		+
		\begin{bmatrix}
			x_{1,m} \\
			x_{2,m}
		\end{bmatrix}, \ t \in [0, 2\pi).
	\end{equation}
	The the geometric parameters $c_m = \sqrt{a_m^2 - b_m^2}$ and $\rho_m = \ln((a_m + b_m)/c_m)$.
\end{definition}

The explicit spectral properties of the NP operator are given by the following theorem.
\begin{theorem}[\cite{AK2016Analysis, chung2014cloaking}]\label{thm:spectral_of_NP}
	The $i$-th order eigenvalues and corresponding eigenfunctions of the NP operator $\mathcal{K}_{w_m}^*$ are given by 
		\begin{gather}
		 \psi_{m,i}^s(t) = \sin(i t)/\Xi(\rho_m, t), \quad
		 \psi_{m,i}^c(t) = \cos(i t)/\Xi(\rho_m, t),\label{eq:eigenfunction_of_Kswm}\\
		\mathcal{K}_{w_m}^*[\psi_{m,i}^s] = -\alpha_{m,i}\psi_{m,i}^s, \quad
		\mathcal{K}_{w_m}^*[\psi_{m,i}^c] = \alpha_{m,i} \psi_{m,i}^c, \label{eq:Kswmpsi}
		\end{gather}
	where $\alpha_{m,i} = e^{-2 i \rho_m}/2$. The single layer $\mathcal{S}_{w_m}$ evaluated on the eigenfunctions yields 
	\begin{equation}\label{eq:Swmpsi}
		\begin{aligned}
			\mathcal{S}_{w_m}[\psi_{m,i}^s] = - \left(0.5 - \alpha_{m,i}\right) \sin i t/i, \
			\mathcal{S}_{w_m}[\psi_{m,i}^c] = - \left(0.5 + \alpha_{m,i}\right) \cos i t/i.
		\end{aligned}
	\end{equation}
	For the special case $i = 0$, we have $\mathcal{S}_{w_m}[\psi_{m,0}^c] = \rho_m + \log(c_m/2)$.
\end{theorem}

Next, we consider the spectral properties of $\mathcal{K}_{w_m}$ (the adjoint of the parameterized NP operator $\mathcal{K}_{w_m}^*$), which has similar results to $\mathcal{K}_{w_m}^*$. The proof can be found in \cref{appendix:spectral_of_Kwm}.
\begin{theorem}\label{thm:spectral_of_Kwm}
	The $i$-th order eigenvalues and corresponding eigenfunctions of parameterized operator $\mathcal{K}_{w_m}$ are given by 
	\begin{gather}
		q_{m,i}^s(t) = \sin(i t)\Xi(\rho_m, t), \quad
		q_{m,i}^c(t) = \cos(i t)\Xi(\rho_m, t),\label{eq:eigenfunction_of_Kwm}\\
		\mathcal{K}_{w_m}[q_{m,i}^s] = -\alpha_{m,i}q_{m,i}^s, \quad 
		\mathcal{K}_{w_m}[q_{m,i}^c] = \alpha_{m,i}q_{m,i}^c, \label{eq:Kwmpsi}
	\end{gather}
	where $\alpha_{m,i} = e^{-2i\rho_m}/2$. Furthermore, we have
	\begin{align}\label{eq:Sswmpsi}
		S^*_{w_m}[\sin(is)] = -(0.5 - \alpha_{m,i})q_{m,i}^s/i, \
		S^*_{w_m}[\cos(is)] = -(0.5 + \alpha_{m,i})q_{m,i}^c/i. 
	\end{align}
	For the special case $i = 0$, we have $S^*_{w_m}[\cos(is)](t) = (\rho_m + \ln(c_m/2))\Xi(\rho_m,t)$.
\end{theorem}

\subsection{RBM for scattering problems}\label{subsec:solve_scattering_prolbem}

In this subsection, we solve the multiple scattering problem \eqref{eq:parameterized_state_equation} using singular splitting strategy combined with the eigenfunctions of the NP operator, which is shape-dependent. For convenience in the subsequent discussion, we relabel the eigenfunctions using a unified index:
\begin{equation*}
	\psi_{m,i} = \psi^s_{m,i}, \ 1 \leq i \leq N/2 - 1 \quad \text{and}\quad \psi_{m,i} = \psi^c_{m,i - N/2}, \ N/2 \leq i \leq N,
\end{equation*}
where $N$ is the total number of retained eigenfunctions (cut-off number). Inspired by the idea of MEM, let $\tilde{\phi}_m^N$ and $\tilde{\varphi}_m^N$ denote the approximations of the densities $\tilde{\phi}_m$ and $\tilde{\varphi}_m$ for the $m$-th particle
\begin{align*}
	\tilde{\phi}_m^N(t) = \sum_{i=1}^{N} c_{m,i}^{\tilde{\phi}} \, \psi_{m,i}(t), \quad
	\tilde{\varphi}_m^N(t) = \sum_{i=1}^{N} c_{m,i}^{\tilde{\varphi}} \, \psi_{m,i}(t).
\end{align*}
 The full density vectors $\tilde{\phi}$ and $\tilde{\varphi}$ are then approximated by the concatenation of their component functions $\tilde{\phi}^N = (\tilde{\phi}_1^N, \ldots, \tilde{\phi}_M^N)^\top$ and $\tilde{\varphi}^N = (\tilde{\varphi}_1^N, \ldots, \tilde{\varphi}_M^N)^\top$, with corresponding coefficiens
$\mathbf{c}^{\tilde{\phi}} = ( \mathbf{c}_1^{\tilde{\phi}}, \ldots, \mathbf{c}_M^{\tilde{\phi}} )^\top$, $\mathbf{c}^{\tilde{\varphi}} = ( \mathbf{c}_1^{\tilde{\varphi}}, \ldots, \mathbf{c}_M^{\tilde{\varphi}} )^\top\in \mathbb{C}^{MN}$, where each $\mathbf{c}_m^{\tilde{\phi}} = (c_{m,1}^{\tilde{\phi}}, \ldots, c_{m,N}^{\tilde{\phi}})^\top \in \mathbb{C}^N$, $m=1,\cdots,M$, and similarly for $\mathbf{c}_m^{\tilde{\varphi}}$. Next, we apply a collocation method to determine unknown coefficients $\mathbf{c}^{\tilde{\phi}}$ and $\mathbf{c}^{\tilde{\varphi}}$. Let $t_j = 2\pi j/N$, for $j = 1, \ldots, N$, be a set of equispaced collocation points on each particle. The fully discretized system of \eqref{eq:parameterized_state_equation} then reads: Find $\mathbf{c}^{\tilde{\phi}}$ and $\mathbf{c}^{\tilde{\varphi}}$ such that
\begin{equation}\label{eq:full_discretization_scattering_problem}
	\begin{cases}
		\mathcal{S}_w^{k_c}[\tilde{\phi}^N](t_j) - \mathcal{S}_w^{k_m}[\tilde{\varphi}^N](t_j) = f_1(t_j), \\[1mm]
		\displaystyle
		\frac{1}{\varepsilon_c} \left( -\frac{1}{2}\mathcal{I} + (\mathcal{K}_w^{k_c})^* \right)[\tilde{\phi}^N](t_j)
		- \frac{1}{\varepsilon_m} \left( \frac{1}{2}\mathcal{I} + (\mathcal{K}_w^{k_m})^* \right)[\tilde{\varphi}^N](t_j) = f_2(t_j).
	\end{cases}
\end{equation}

To address the singularity, the kernel of the boundary operator can be decomposed as  $G^k(z)=G(z)+\widehat{G}^k(z)$, where $G(z)$ represents the singular part and $\widehat{G}^k(z)$ is the non-singular part which depends on the wave number $k$. Using this property, the boundary operator can be decomposed into singular and non-singular parts. The non-singular parts are defined as $\widehat{\mathcal{S}}_{w_m}^k = \mathcal{S}_{w_m}^k - \mathcal{S}_{w_m}$ and $(\widehat{\mathcal{K}}_{w_m}^k)^* =  (\mathcal{K}_{w_m}^{k})^* - \mathcal{K}_{w_m}^*$. As shown in \eqref{eq:limit_of_hatGk}, the kernel of  $\widehat{\mathcal{S}}_{w_m}^k$ and $(\widehat{\mathcal{K}}_{w_m}^k)^*$ remain non-singular as $t \to s$.

Let $\mathbf{S}^k_w$ denote the discrete form of the single-layer potential operator evaluated on eigenfunctions, which can be decomposed into a singular part and a non-singular part, i.e. $\mathbf{S}^k_w = \mathbf{S}_w + \widehat{\mathbf{S}}^k_w$, with
\begin{equation*}
	\begin{aligned}
		& \mathbf{S}_w= \diag[\mathcal{S}_{w_1}[\{\psi_{1,i}\}](\{t_j\}), \cdots, \mathcal{S}_{w_M}[\{\psi_{M,i}\}](\{t_j\})]\in \mathbb{C}^{MN \times MN},\\
		&\widehat{\mathbf{S}}^k_w = \begin{bmatrix}
			\widehat{\mathcal{S}}^k_{w_1}[\{\psi_{1,i}\}](\{t_j\}) & \cdots & \mathcal{S}^k_{w_1,w_M}[\{\psi_{M,i}\}](\{t_j\}) \\
			\vdots & \ddots & \vdots \\
			\mathcal{S}^k_{w_M,w_1}[\{\psi_{1,i}\}](\{t_j\}) & \cdots & \widehat{\mathcal{S}}^k_{w_M}[\{\psi_{M,i}\}](\{t_j\}) \\
		\end{bmatrix}\in \mathbb{C}^{MN \times MN}.
	\end{aligned}
\end{equation*}
Here, the block matrix in $\mathbf{S}_w$ is 
\begin{equation*}
	\begin{aligned}
		\mathcal{S}_{w_m}[\{\psi_{m,i}\}](\{t_j\}) = \begin{bmatrix}
		\mathcal{S}_{w_m}[\psi_{m,1}](t_1) & \cdots & \mathcal{S}_{w_m}[\psi_{m,N}](t_1) \\
		\vdots & \ddots & \vdots \\
		\mathcal{S}_{w_m}[\psi_{m,1}](t_N) & \cdots & \mathcal{S}_{w_m}[\psi_{m,N}](t_N) \\
		\end{bmatrix} \in \mathbb{C}^{N \times N}.
	\end{aligned}
\end{equation*}
Note that each term in singular part can be computed explicitly using formulas presented in \Cref{thm:spectral_of_NP}. On other hand, each term in the non-singular part can be computed via the trapezoidal rule with high accuracy, since the kernel is smooth \cite{kanwal2013linear}. Similarly, $(\mathbf{K}^k_w)^*$ can also be decomposed into a singular and a non-singular part. Thus, the fully discretized form \eqref{eq:full_discretization_scattering_problem} can be reformulated as the linear system:
\begin{equation}\label{eq:linear_system_of_scattering_problem}
	\mathbf{M}[w,\lambda]\, \mathbf{c} = \mathbf{f}[w,\lambda],
\end{equation}
where $\mathbf{c}=[\mathbf{c}^{\tilde{\phi}}, \mathbf{c}^{\tilde{\varphi}}]^{\top} $ and $\mathbf{f}=[f_1(\{t_j\}), f_2(\{t_j\})]^{\top} $. Moreover, $\mathbf{M}[w,\lambda] = \mathbf{M}_1[w,\lambda] + \mathbf{M}_2[w,\lambda]\in \mathbb{C}^{2MN \times 2MN}$ is composed of a singular part and a non-singular part. Specifically, they are given by:
\begin{equation}\label{eq:decomposition_of_scattering_matrix}
\resizebox{0.9\textwidth}{!}{$
	\mathbf{M}_1 =
	\begin{bmatrix}
		\mathbf{S}_w & -\mathbf{S}_w \\
		(-\frac{1}{2} \mathbf{I} + \mathbf{K}_w^*)/\varepsilon_c &
		 -(\frac{1}{2} \mathbf{I} + \mathbf{K}_w^*)/\varepsilon_m
	\end{bmatrix}
	,\
	\mathbf{M}_2 =
	\begin{bmatrix}
		\widehat{\mathbf{S}}^{k_c}_w & -\widehat{\mathbf{S}}^{k_m}_w \\
		(\widehat{\mathbf{K}}^{k_c}_w)^*/\varepsilon_c &
		-(\widehat{\mathbf{K}}^{k_m}_w)^*/\varepsilon_m
	\end{bmatrix}.
    $}
\end{equation}

\begin{remark}
	For broadband computations, the block matrices of the singular part $\mathbf{S}_w$ and $(\mathbf{K}_w)^*$ can be precomputed. This precomputation allows for the non-singular parts to be evaluated separately for each wavelength, significantly reducing the overall computational cost.
\end{remark}

\subsection{RBM for adjoint problem}\label{subsec:solve_adjoint_equation}

In this subsection, we consider the numerical solution of the adjoint equation \eqref{eq:adjoint_equation}. The approach follows the same strategy used for solving the forward scattering problem \eqref{eq:parameterized_state_equation}, utilizing shape-dependent eigenfunctions for discretization. To address the singularity of the boundary integral operators, we also decompose $(\mathcal{S}_{w_m}^k)^*$ and $\mathcal{K}_{w_m}^k$ into singular and non-singular components, where the non-singular parts are defined as: $(\widehat{\mathcal{S}}_{w_m}^k)^*= (\mathcal{S}_{w_m}^k)^* - \mathcal{S}^*_{w_m}$ and $\widehat{\mathcal{K}}_{w_m}^k =\mathcal{K}_{w_m}^{k} - \mathcal{K}_{w_m}$.

Next, we utilize the properties of $\mathcal{S}_{w_m}^{*}$ and $\mathcal{K}_{w_m}$ given in \Cref{thm:spectral_of_Kwm} to solve the adjoint equation \eqref{eq:adjoint_equation}. For convenience, we define the basis functions as follows:
\begin{equation*}
		\psi^{\tilde{p}}_{m,i} = \begin{cases}
			\begin{aligned}
				&\sin(it),&& 1\leq i\leq N/2-1,\\
				&\cos((i-N/2)t), && N/2\leq i\leq N,
			\end{aligned}
		\end{cases}
		\psi^{\tilde{q}}_{m,i} = \begin{cases}
			\begin{aligned}
				&q_{m,i}^s,&& 1\leq i\leq N/2-1,\\
				&q^c_{m,i - N/2}, && N/2\leq i\leq N.
			\end{aligned}
		\end{cases}
\end{equation*}
Let $\tilde{p}^N_m$, $\tilde{q}^N_m$ denote the approximations of $\tilde{p}_m$ and $\tilde{q}_m$, respectively, 
\begin{equation*}
	\tilde{p}^N_m(t)= \sum_{n=1}^{N}d^{\tilde{p}}_{m,n} \psi^{\tilde{p}}_{m,n}(t),\quad
	\tilde{q}^N_m(t)= \sum_{n=1}^{N}d^{\tilde{q}}_{m,n} \psi^{\tilde{q}}_{m,n}(t).
\end{equation*}
Then $\tilde{p}^N = (\tilde{p}^N_1, \ldots, \tilde{p}^N_M)^\top$ and $\tilde{q}^N = (\tilde{q}^N_1, \ldots, \tilde{q}^N_M)^\top$ are the approximations of $\tilde{p}$ and $\tilde{q}$. Thus, the fully discrete form reads: find $\mathbf{d}^{\tilde{p}}$ and $\mathbf{d}^{\tilde{q}}$ such that
\begin{equation}\label{eq:full_discretization_adjoint_equation}
	\begin{cases}
		\displaystyle (\mathcal{S}_w^{k_c})^*[\tilde{p}^N](t_j) + \frac{1}{\overline{\varepsilon_c}}\left(-\frac{1}{2} I+\mathcal{K}_w^{k_c}\right)[\tilde{q}^N](t_j) = g_1(t_j), \\[2mm]
		\displaystyle -(\mathcal{S}_w^{k_m})^*[\tilde{p}^N](t_j) - \frac{1}{\varepsilon_m}\left(\frac{1}{2} I+\mathcal{K}_w^{k_m}\right)[\tilde{q}^N](t_j) = g_2(t_j).
	\end{cases}
\end{equation}

Following similar arguments as RBM framework for solving scattering problem, the full discretization \eqref{eq:full_discretization_adjoint_equation} can be reformulated as the following linear system:
\begin{align}\label{eq:linear_system_of_adjoint_equation}
	\mathbf{T}[w,\lambda]\,\mathbf{d} = \mathbf{g}[w,\lambda],
\end{align}
where the unknown vector and right hand side are given by $\mathbf{d} = [\mathbf{d}^{\tilde{p}}, \mathbf{d}^{\tilde{q}}]^{\top}$, $\mathbf{g} = [g_1(\{t_j\}), g_2(\{t_j\})]^{\top}$, and $\mathbf{T}[w,\lambda] = \mathbf{T}_1[w,\lambda] + \mathbf{T}_2[w,\lambda]\in \mathbb{C}^{2MN \times 2MN}$, with
\begin{equation}\label{eq:decomposition_of_adjoint_matrix}
\resizebox{0.9\textwidth}{!}{$
	\begin{aligned}
		\mathbf{T}_1[w,\lambda] =\begin{bmatrix}
			(\mathbf{S}_w)^* &  (-\frac{1}{2} \mathbf{I} + \mathbf{K}_w)/\overline{\varepsilon_c} \\[2mm]
			-(\mathbf{S}_w)^*  & -(\frac{1}{2} \mathbf{I} + \mathbf{K}_w)/\varepsilon_m
		\end{bmatrix} , \
		\mathbf{T}_2[w,\lambda] =\begin{bmatrix}
			(\widehat{\mathbf{S}}^{k_c}_w)^*  &   \widehat{\mathbf{K}}^{k_c}_w/\overline{\varepsilon_c} \\[2mm]
			-(\widehat{\mathbf{S}}^{k_m}_w)^*  &  -\widehat{\mathbf{K}}^{k_m}_w/\varepsilon_m
		\end{bmatrix} .
	\end{aligned}
    $}
\end{equation}

\subsection{Derivative of boundary operator}\label{subsec:derivative_of_BO}

In this subsection, we consider the derivative of the boundary operators with respect to the parameters \(w\), as outlined in \eqref{eq:dEdw_and_dFdw}, and address the singularity using the eigenfunctions of the NP operator.

According to the definition in \eqref{eq:derivative_of_boundary_op}, the derivatives of the non-diagonal terms  with respect to $w_n$ or $w_m$ are non-singular. These derivatives can be evaluated numerically using the trapezoidal rule with high accuracy. Next, we focus on the derivative of the diagonal terms. By applying singular splitting again, we have
\begin{equation*}
	\begin{aligned}
		\frac{\partial \mathcal{S}_{w_m}^{k}}{\partial w_m} [\tilde{\psi}_m]=  \frac{\partial \mathcal{S}_{w_m}}{\partial w_m}[\tilde{\psi}_m]+\frac{\partial \widehat{\mathcal{S}}_{w_m}^{k}}{\partial w_m}[\tilde{\psi}_m],\
		\frac{\partial (\mathcal{K}_{w_m}^{k})^*}{\partial w_m}[\tilde{\psi}_m] =  \frac{\partial \mathcal{K}_{w_m}^*}{\partial w_m}[\tilde{\psi}_m] +\frac{\partial (\widehat{\mathcal{K}}_{w_m}^{k})^*}{\partial w_m}[\tilde{\psi}_m].
	\end{aligned}
\end{equation*}

First, we consider the derivative of the singular operator:  
\begin{equation}\label{eq:singular_part_of_BDO_derivative}
	\begin{aligned}
	&\frac{\partial \mathcal{S}_{w_m}}{\partial w_m}[\tilde{\psi}_m]=  \sum_{i=1}^{N} c_{m,i} \frac{\partial \mathcal{S}_{w_m}}{\partial w_m}[\psi_{m,i}],\quad
	\frac{\partial \mathcal{K}_{w_m}^*}{\partial w_m}[\tilde{\psi}_m] =   \sum_{i=1}^{N}c_{m,i} \frac{\partial \mathcal{K}_{w_m}^*}{\partial w_m}[\psi_{m,i}],
	\end{aligned}
\end{equation}
where we assume $\tilde{\psi}_m = \sum_{i=1}^{N} c_{m,i} \psi_{m,i}$. By using the following identities,  
\begin{align*}
	\frac{\partial \mathcal{S}_{w_m}}{\partial w_m}[\psi_{m,i}] =  \frac{\partial \mathcal{S}_{w_m}[\psi_{m,i}]}{\partial w_m}-\mathcal{S}_{w_m}[\frac{\partial \psi_{m,i}}{\partial w_m}],\
	\frac{\partial \mathcal{K}_{w_m}^*}{\partial w_m}[\psi_{m,i}] =  \frac{ \partial \mathcal{K}_{w_m}^*[\psi_{m,i}]}{\partial w_m}- \mathcal{K}_{w_m}^*[\frac{\partial \psi_{m,i}}{\partial w_m}],
\end{align*}
the terms on the right-hand side can be calculated explicitly, as shown in \Cref{sec:appendix_Derivative_of_BO}.

Second, we consider the derivative of the non-singular operator:
\begin{equation*}
	\begin{aligned}
		&\frac{\partial \widehat{\mathcal{S}}_{w_m}^{k}}{\partial w_m}[\tilde{\psi}_m] 
		= \int_{0}^{2\pi} \left(\frac{\partial \widehat{G}^{k}(z(t,s;w_m))}{\partial w_m}|y^{\prime}(s;w_m)| +\widehat{G}^{k}(z(t,s;w_m)) \frac{\partial |y^{\prime}(s;w_m)|}{\partial w_m} \right) \tilde{\psi}_m(s)\mathrm{d}s,\\
		&\frac{\partial (\widehat{\mathcal{K}}_{w_m}^{k})^*}{\partial w_m}[\tilde{\psi}_m] =  \int_{0}^{2\pi} \frac{\partial }{\partial w_m}\left( \frac{\partial \widehat{G}^{k}(z(t,s;w_m))}{\partial \nu_x(t;w_m)} |y^{\prime}(s;w_m)| \right) \tilde{\psi}_m(s)\mathrm{d}s \\
		&=  \int_{0}^{2\pi} \left[\frac{\partial}{\partial w_m} \left(\frac{\partial \widehat{G}^{k}(z(t,s;w_m))}{\partial \nu_x(t;w_m)}\right)|y^{\prime}(s;w_m)|+\frac{\partial \widehat{G}^{k}(z(t,s;w_m))}{\partial \nu_x(t;w_m)} \frac{\partial |y^{\prime}(s;w_m)|}{\partial w_m} \right]  \tilde{\psi}_m(s)\mathrm{d}s. 
	\end{aligned}
\end{equation*}
As shown in \Cref{appendix:asymptotic_of_Green}, the kernels appearing in the integrals are non-singular. Hence, they can also be approximated  using the trapezoidal rule.

\section{Numerical Algorithm}\label{sec:algorithm}

In this section, we present the numerical algorithm for solving the optimal design problem \eqref{eq:parameterized_optimal_design}. We begin with the formulation of a strategy for constructing an initial guess, which is critical for ensuring the the convergence to a good solution. Subsequently, a comprehensive description of the proposed algorithm is provided, which is based on the gradient descent approach.

\subsection{Initial guess generation}

This subsection introduces a physics-informed methodology to generate an initial guess for the optimal design problem described in \eqref{eq:parameterized_optimal_design}. The approach begins with the definition of an offline-constructed dataset: $\{A(w_{\ell}, \lambda)\}_{\ell=1}^L$, where $w_{\ell} = (a_{\ell}, b_{\ell}, \theta_{\ell}, 0, 0)^\top$ and $\lambda \in \Lambda$. Here, $A(w_{\ell}, \lambda)$ denotes the absorptance of a single ellipse centered at the origin. The parameters are subject to the following constraints: $a_{\min} \leq a_{\ell} \leq a_{\max}$, $\eta_{\min} a_{\ell} \leq b_{\ell} \leq \eta_{\max} a_{\ell}$, and $0 \leq \theta_{\ell} \leq \pi/2$. Specifically, $a_{\min}$ and $a_{\max}$ define the permissible range for the semi-major axis length, $\eta_{\min}$ and $\eta_{\max}$ establish the aspect ratio bounds, and $\theta_{\ell}$ is restricted to $[0, \pi/2]$ due to symmetry considerations.
\begin{remark}
	The upper and lower bounds on the semi-major axis parameter $a$, namely $a_{\max}$ and $a_{\min}$ specify the range of particle sizes. Additionally, the bounds on the aspect ratio parameter $\eta$, given by $\eta_{\min}$ and $\eta_{\max}$, are imposed to prevent the ellipse from becoming too singular.
\end{remark}

In the subwavelength regime \cite{bohren2008absorption} and under the weakly interacting assumption based on multiple scattering theory \cite{martin2006multiple}, the total absorptance can be approximated as the superposition of contributions from individual particles. Given a target absorptance $A^{\mathrm{tar}}(\lambda)$, the optimal particle configuration is determined by solving the following quadratic integer programming problem:
\begin{align}\label{eq:matrix_integer_opt}
	c^* = \mathop{\arg\min}_{c \in \mathbb{Z}_+^L} \| D(\lambda) c - A^{\mathrm{tar}}(\lambda) \|^2_{L^2(\Lambda)},
\end{align}
where the dataset $D(\lambda) = [A(w_1, \lambda), \ldots, A(w_L, \lambda)]$ and $c = [c_1, \ldots, c_L]^\top$. Here, $c_{\ell}$ represents the multiplicity of particles characterized by the parameters $w_{\ell}$.  The primary focus is on obtaining a computationally feasible initial guess rather than an exact optimal solution since solving \eqref{eq:matrix_integer_opt} exactly is NP-hard.

First, the integer constraints in \eqref{eq:matrix_integer_opt} are relaxed to formulate a continuous form:
\begin{equation}\label{eq:relax_form_of_matrix_integer_opt}
	c^\dagger = \mathop{\arg\min}_{c \in \mathbb{R}^L_+} \| D(\lambda) c - A^{\mathrm{tar}}(\lambda) \|^2_{L^2(\Lambda)}.
\end{equation}
The relaxed problem is solved using quadratic programming, yielding the solution $c^\dagger$. Subsequently, this solution is rounded element-wise to the nearest integers to obtain $c^\ddagger = \lfloor c^\dagger \rceil$. Second, the initial approximation $c^\ddagger$ is further refined using a heuristic algorithm to produce the final solution, denoted as $c^*$. Third, the initial particle configuration is generated based on $c^*$ as follows:
\begin{itemize}[leftmargin=*]
	\item The total number of particles is calculated as: $M = \sum_{\ell=1}^L c_{\ell}^*$.
	
	\item Geometric parameters are assigned by:
	\begin{center}
		$a^{(0)}_{m} = a_{\ell}$, $b^{(0)}_{m} = b_{\ell}$, $\theta^{(0)}_{m} = \theta_{\ell}$ , for $m \in (\sum_{i=0}^{\ell-1} c_i, \sum_{i=1}^\ell c_i]$, $1 \leq \ell \leq L$.
	\end{center}
	
	\item Particle positions are distributed on a uniform grid:
	\begin{center}
		$x^{(0)}_{1,m} = ( i - (1 + N_{x_1})/2 ) \Delta_1$, $x^{(0)}_{2,m} = ( j - (1 + N_{x_2})/2 ) \Delta_2$, $m = (j - 1) N_{x_1} + i$,
	\end{center}
	where $\Delta_1$ and $\Delta_2$ are the grid spacings, and $N_{x_1} N_{x_2} \leq M$ specifies the grid size.
	
	\item The complete initial guess is represented as:
	\begin{align}\label{eq:w0}
		w^{(0)} = (w^{(0)}_1, \ldots, w^{(0)}_M), \quad
		w^{(0)}_m = (a^{(0)}_{m}, b^{(0)}_{m}, \theta^{(0)}_{m}, x^{(0)}_{1,m}, x^{(0)}_{2,m}).
	\end{align}
\end{itemize}
The initialization strategy is summarized in \Cref{algo:initial_guess}.
\begin{algorithm}[htbp]
	\caption{Initial Guess Generation}
	\label{algo:initial_guess}
	\begin{algorithmic}[1]
		\REQUIRE Target absorptance $A^{\mathrm{tar}}(\lambda)$ and dataset $\{A(w_{\ell}, \lambda)\}_{\ell=1}^L$.
		
		\STATE Solve the relaxed problem \eqref{eq:relax_form_of_matrix_integer_opt} to obtain $c^\dagger$ and $c^\ddagger \gets \lfloor c^\dagger \rceil$.
		
		\STATE Refine $c^\ddagger$ using a heuristic algorithm to obtain $c^*$.
		
		\STATE Construct the initial configuration $w^{(0)}$ using $c^*$ via \eqref{eq:w0}.
		
		\ENSURE Initial particle configuration $w^{(0)}$.
	\end{algorithmic}
\end{algorithm}

\subsection{Optimal design}

In this subsection, we employ the gradient descent method, utilizing the initial guess from \cref{algo:initial_guess}, to solve the optimal design problem \eqref{eq:parameterized_optimal_design}. Each iteration's computation is decomposed into a wavelength-independent singular part and a wavelength-dependent non-singular part. This decomposition enables efficient computation of the scattering problem, the adjoint problem, and the derivatives of boundary operators. The procedure consists of the following steps:

\textbf{Step 1: Singular Components.}  
For fixed shape parameters $w = w^{(i)}$, the primary wavelength-independent singular components $\mathbf{S}_w$ and $\mathbf{K}_w^*$ (for the scattering problem), as well as $\mathbf{S}_w^*$ and $\mathbf{K}_w$ (for the adjoint problem), are precomputed. Note that the singular matrices $\mathbf{M}_1[w, \lambda]$ and $\mathbf{T}_1[w, \lambda]$ depend on the wavelength-dependent material parameters $\varepsilon_m$ and $\varepsilon_c$, and can be computed easily via equations~\eqref{eq:decomposition_of_scattering_matrix} and~\eqref{eq:decomposition_of_adjoint_matrix}. Moreover, the singular part of the derivative of the boundary operator acting on its eigenfunction, i.e., $\frac{\partial \mathcal{S}_{w_m}}{\partial w_m}[\psi_{m,i}]$ and $\frac{\partial \mathcal{K}_{w_m}^*}{\partial w_m}[\psi_{m,i}]$, can be precomputed. Note that, for each solution obtained by the RBM, the singular parts of $\mathcal{E}_w[\tilde{\phi}, \tilde{\varphi}, w]$ and $\mathcal{F}_w[\tilde{\phi}, \tilde{\varphi}, w, \lambda]$ can be computed easily using equation~\eqref{eq:singular_part_of_BDO_derivative}.

\textbf{Step 2: Non-Singular Components.}  
For fixed $w = w^{(i)}$ and each $\lambda = \lambda_j$, we compute the non-singular parts as follows. First, the non-singular matrices $\mathbf{M}_2[w, \lambda]$ and $\mathbf{T}_2[w, \lambda]$ are determined via the trapezoidal rule. Note that the non-singular operators for the scattering problem, which have kernels $\widehat{G}^k(z(t, s; w_m))$ and $G^k(z(t, s; w_n, w_m))$, are calculated for $k \in \{k_m, k_c\}$. Using the symmetry property $G^{-k}(z) = \overline{G^k(z)}$, the kernels $\widehat{G}^{-k}(z(t, s; w_m))$ and $G^{-k}(z(t, s; w_n, w_m))$ for the adjoint problem can be obtained directly. Following this, the discretized scattering system \eqref{eq:full_discretization_scattering_problem} and the adjoint system \eqref{eq:full_discretization_adjoint_equation} are solved to obtain the necessary fields. Finally, the gradient components $A_w(\tilde{\varphi}, w, \lambda)$ are calculated using equation~\eqref{eq:dAdw}, while the non-singular parts of $\mathcal{E}_w[\tilde{\phi}, \tilde{\varphi}, w]$ and $\mathcal{F}_w[\tilde{\phi}, \tilde{\varphi}, w, \lambda]$ are also computed using the trapezoidal rule.

\textbf{Step 3: Update Parameters.}  
First, we compute the full gradient $\mathcal{J}^\prime(w^{(i)})$ over the broadband via equation~\eqref{eq:full_gradient}. Second, we update the shape parameters using gradient descent:
\begin{equation}\label{eq:gradient_descent}
	\widetilde{w}^{(i+1)} = w^{(i)} - \beta \mathcal{J}^\prime(w^{(i)}).
\end{equation}  
Third, we project the updated parameters onto the constraint set:
\begin{equation}\label{eq:projection}
	\begin{aligned}
		a^{(i+1)}_m &= \mathcal{P}_{[a_{\min}, a_{\max}]}(\widetilde{a}^{(i+1)}_m), \quad
		\theta^{(i+1)}_m = \mathcal{P}_{[0, 2\pi]}(\widetilde{\theta}^{(i+1)}_m), \\
		\eta^{(i+1)}_m &= \mathcal{P}_{[\eta_{\min}, \eta_{\max}]}(\widetilde{\eta}^{(i+1)}_m), \quad
		b^{(i+1)}_m = \eta^{(i+1)}_m a^{(i+1)}_m.
	\end{aligned}
\end{equation}  
where the projection operator $\mathcal{P}_{[c, d]}(x) = \max(c,\min(x,d))$. Then, the updated particle parameters are 
\begin{equation}\label{eq:updated_parameters}
	w^{(i+1)}_m = (a^{(i+1)}_m, b^{(i+1)}_m, \theta^{(i+1)}_m, x^{(i+1)}_{1,m}, x^{(i+1)}_{2,m}).
\end{equation}

The complete procedure for the optimal design of broadband absorbers using multiple particles is summarized in \cref{algo:optimal_design}.

\begin{algorithm}[ht]
	\caption{Optimal Design of Broadband Absorber via Multiple Particles}
	\label{algo:optimal_design}
	\begin{algorithmic}[1]
		\REQUIRE Target absorptance $A^{\mathrm{tar}}(\lambda)$; 
		off-line absorptance dataset $\{ A(w_{\ell}, \lambda) \}_{\ell=1}^L$; 
		number of basis functions per particle $N$; 
		number of wavelength quadrature points $N_{\lambda}$; 
		maximum iterations $N_{\mathrm{iter}}$; 
		bounds $a_{\min}$, $a_{\max}$ for the semi-major axis; 
		aspect ratio bounds $\eta_{\min}$, $\eta_{\max}$;
		step size $\beta$.
		
		\STATE Generate an initial guess $w^{(0)}$ via \Cref{algo:initial_guess}.
		
		\FOR{$i = 0 : N_{\mathrm{iter}}$}
		\STATE Precompute the following wavelength-independent singular components: 
		\begin{center}
			Scattering: $\mathbf{S}_w$, $\mathbf{K}_w^*$; 
			Adjoint: $\mathbf{S}_w^*$, $\mathbf{K}_w$; 
			Derivative: $\frac{\partial \mathcal{S}_{w_m}}{\partial w_m}[\psi_{m,i}]$, $\frac{\partial \mathcal{K}_{w_m}^*}{\partial w_m}[\psi_{m,i}]$.
		\end{center}
		
		\FOR{$j = 0 : N_{\lambda}$}
		\STATE Compute the singular parts $\mathbf{M}_1[w^{(i)},\lambda_j]$ and $\mathbf{T}_1[w^{(i)},\lambda_j]$ via~\eqref{eq:decomposition_of_scattering_matrix},~\eqref{eq:decomposition_of_adjoint_matrix}.
		
		\STATE Compute non-singular parts  $\mathbf{M}_2[w^{(i)},\lambda_j]$ and $\mathbf{T}_2[w^{(i)},\lambda_j]$ via \eqref{eq:decomposition_of_scattering_matrix}, \eqref{eq:decomposition_of_adjoint_matrix}.
		
		\STATE Solve the scattering problem \eqref{eq:parameterized_state_equation} to obtain solutions $\tilde{\phi}^N$ and $\tilde{\varphi}^N$.
		
		\STATE Solve the adjoint problem \eqref{eq:adjoint_equation} to obtain solutions $\tilde{p}^N$ and $\tilde{q}^N$.
		
		\STATE Compute the derivatives of these operators via equations~\eqref{eq:dAdw}-\eqref{eq:dEdw_and_dFdw}:
		\begin{center}
			$A_w(\tilde{\varphi}^N, w^{(i)}, \lambda_j)$, $\mathcal{E}_w[\tilde{\phi}^N, \tilde{\varphi}^N, w^{(i)}, \lambda_j]$, $\mathcal{F}_w[\tilde{\phi}^N, \tilde{\varphi}^N, w^{(i)}, \lambda_j]$.
		\end{center}
		
		\ENDFOR
		
		\STATE Compute the full gradient $\mathcal{J}^\prime(w^{(i)})$ via ~\eqref{eq:full_gradient}.
		\STATE Perform gradient descent to obtain $\widetilde{w}^{(i+1)}$ via ~\eqref{eq:gradient_descent}.
		\STATE Project $\widetilde{w}^{(i+1)}$ onto the constraints to obtain $w^{(i+1)}_m$ via \eqref{eq:projection}.
		\ENDFOR
		\ENSURE Final optimized design $w^{(N_{\mathrm{iter}})}$.
	\end{algorithmic}
\end{algorithm}

\section{Numerical Experiments}\label{sec:numerical}

In this setting, we assume that the homogeneous surrounding medium has a constant relative electric permittivity $\varepsilon_{r,m} = 1$ and a constant magnetic permeability $\mu_{r,m} = 1$, both independent of the wavelength. Additionally, the magnetic permeability of the nanoparticles is set to $\mu_{r,c} = 1$. In contrast, the relative electric permittivity $\varepsilon_{r,c}$ is given by $\varepsilon_{r,c}(\omega) = 1 - \omega_p^2/\omega(\omega + i\tau)$ (Drude model \cite{sarid2010modern}), where $\omega_p$ is the plasma frequency of the bulk material, and $\tau > 0$ denotes the damping coefficient. In our simulations, we consider nanoparticles made of silver, with material parameters $\omega_p = 7.613$~eV and $\tau = 0.048$~eV.

In the following numerical examples, unless otherwise specified, we assume the incident angle $\theta_0=0$, i.e. $d = (1, 0)$. The received energy is measured on a partial circular arc characterized by $\overline{\theta} = \theta_0 = 0$, $\Delta \theta = \pi/4$, and radius $R = 1500$~nm. The number of eigenfunctions used for each particle is set to $N = 10$. The inter-particle spacings are chosen to be $\Delta_1 = \Delta_2 = 80$~nm.

\subsection{Absorptance vs shape, rotation, and multi-particle effects}\label{example:shape_angle_combination}
	
	In this example, we investigate the influence of shape, rotation angle, and the presence of multiple particles on absorptance.

\begin{figure}[htbp]
		\centering
		\includegraphics[width=.95\textwidth]{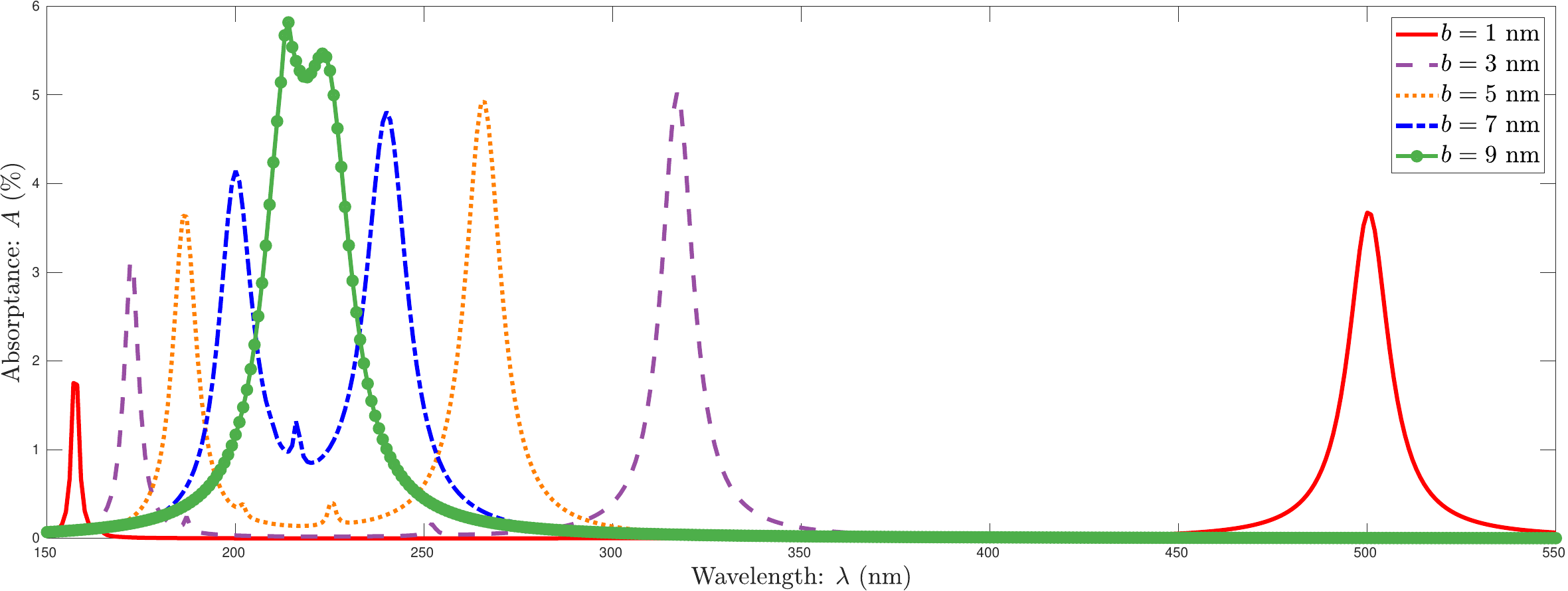}
		\caption{Absorptance versus semi-minor axis $b$ for an origin-centered elliptical nanoparticle ($a = 10$ nm, $\theta = \pi/4$). As $b$ decreases, the two resonance peaks move further apart and their magnitudes diminish.}
		\label{fig:ratio}
	\end{figure}

	First, we examine how absorptance depends on particle shape by fixing the semi-major axis at $a = 10$~nm and varying the semi-minor axis $b$. A single particle is placed at the origin and rotated by an angle $\theta = \pi/4$. As shown in \cref{fig:ratio}, decreasing the semi-minor axis $b$—which produces a flatter elliptical shape—causes a single resonance frequency to split into two distinct resonance frequencies. As $b$ is further decreased, the two resonances become more widely separated, and the peak magnitudes are reduced accordingly.
	
\begin{figure}[htbp]
		\centering
		\includegraphics[width=.95\textwidth]{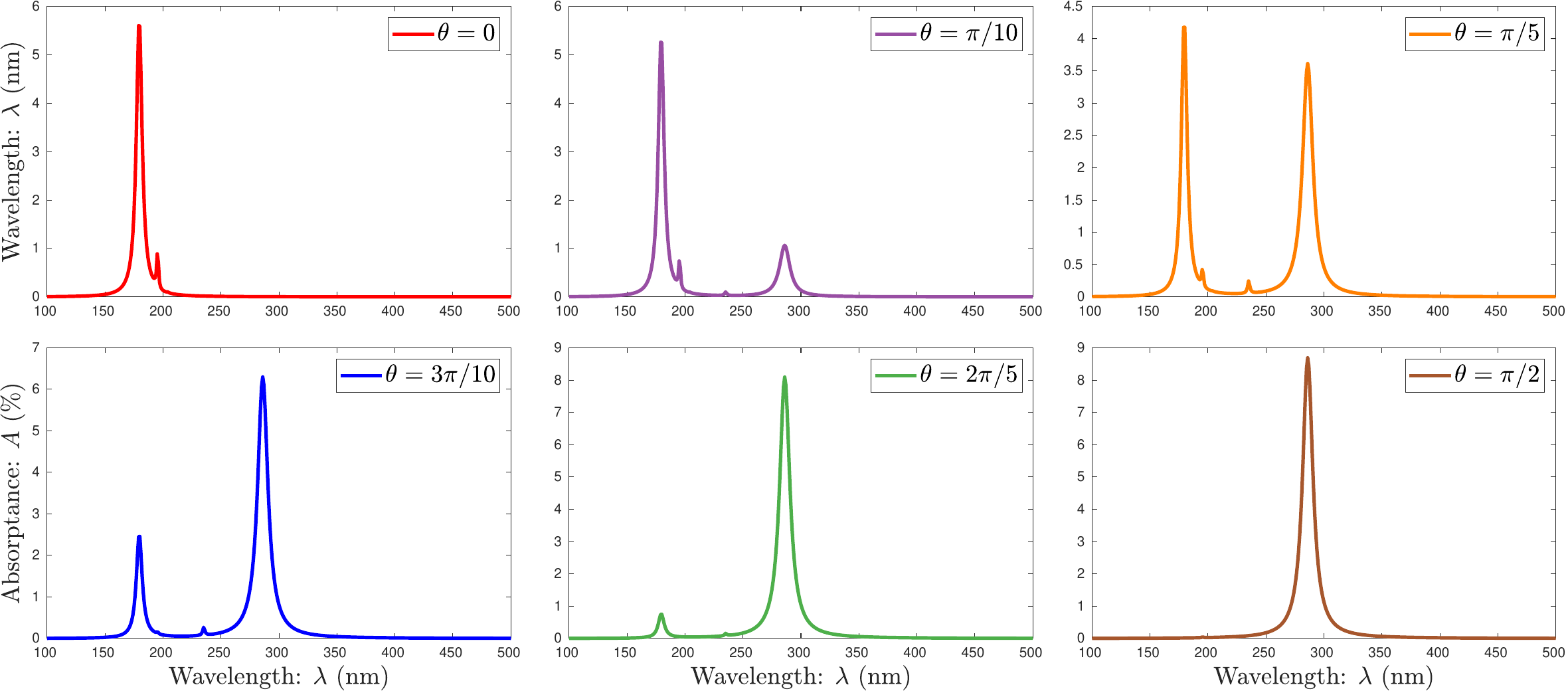}
		\caption{Absorptance versus rotation angle $\theta$ for an origin-centered elliptical nanoparticle ($a = 10$~nm, $b = 4$~nm). This indicates that varying the rotation angle does not shift the resonance frequencies, but it does affect the magnitudes of absorptance at these frequencies.}
		\label{fig:angle}
	\end{figure}

	Second, we examine how absorptance depends on the rotation angle $\theta$, with a single particle fixed at the origin. The numerical results shown in \cref{fig:angle} demonstrate that varying the rotation angle does not shift the resonance frequencies, but it does affect the magnitudes of absorptance at these frequencies.
	
\begin{figure}[htbp]
		\centering
		\includegraphics[width=0.95\textwidth]{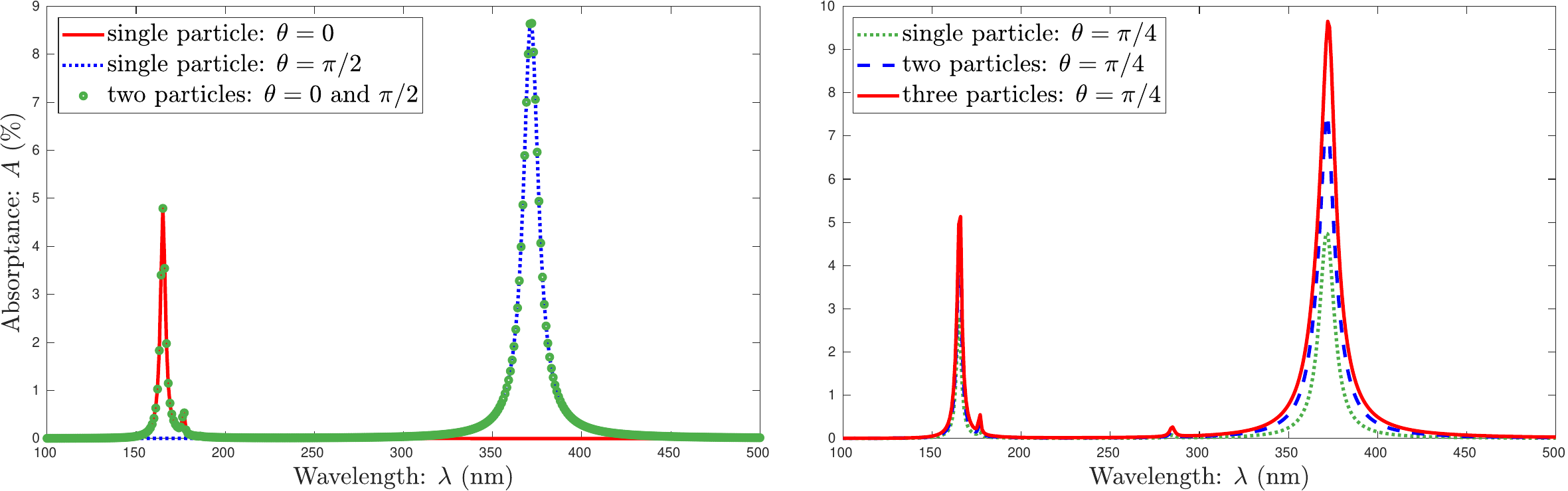}
		\caption{Absorptance for a single particle and combinations of multiple particles ($a = 10$~nm, $b = 2$~nm). Left: Different particles—total absorptance closely approximates the sum of the absorptance of individual particles. Right: Identical particles—total absorptance does not equal the sum of the absorptance of individual particles.}
		\label{fig:weak_scattering}
	\end{figure}

	Third, we examine the weak scattering effect in configurations involving multiple particles, as illustrated in \cref{fig:weak_scattering}. In the case of \emph{different particles}, where two particles have the same shape but different rotation angles, the total absorptance closely approximates the sum of the absorptance of each particle, as shown in the left panel of \cref{fig:weak_scattering}. In contrast, for \emph{identical particles}—where multiple particles share the same shape and orientation—the right panel of \cref{fig:weak_scattering} shows that the total absorptance does not equal the sum of the absorptance of the individual particles. Instead, interactions between the particles enhance the absorptance magnitude at certain resonance frequencies.
	
\subsection{Efficiency of RBM}\label{example:effectiveness_for_RBM}

	In this example, we evaluate the efficiency of our proposed RBM for solving both the scattering and adjoint problems in the cases of a single particle and multiple particles. For comparison, we use the classical Nystr\"{o}m method (cf.~\cite{CK2019Inverse,kanwal2013linear}) as the reference solution, with the number of discretization points set to $200$ for each particle.
	
\begin{figure}[htbp]
		\centering
		\includegraphics[width=.95\textwidth]{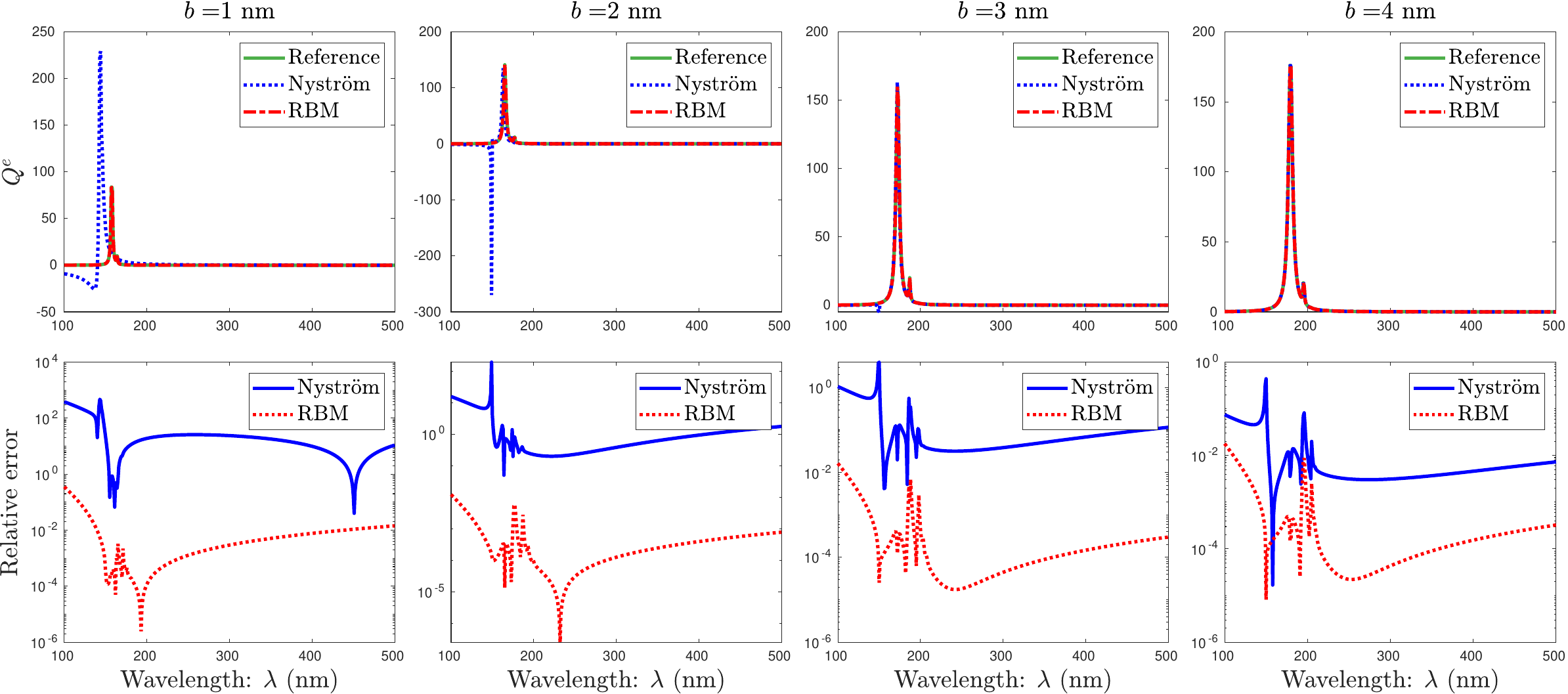}
		\caption{ Extinction cross section $Q^e$ and its relative error $|Q^e_{\mathrm{ref}}-Q^e|/|Q^e_{\mathrm{ref}}|$ for versus semi-minor axis $b$ with an origin-centered elliptical nanoparticle ($a=10$ nm, $\theta=0$, and the number of eigenfuncions $N=10$).}
		\label{fig:scattering_single}
	\end{figure}

	First, we show the efficiency of RBM for solving the scattering problem and computing the extinction cross section $Q^e$, as defined in~\eqref{eq:cross_section}. For the single-particle case, the results in \cref{fig:scattering_single} demonstrate that RBM achieves relative errors several orders of magnitude smaller than those of the Nystr\"{o}m method. Moreover, as the semi-minor axis $b$ increases, the error associated with the Nystr\"{o}m method decreases; however, its performance deteriorates in the flat-particle regime where $a \gg b$.

\begin{figure}[htbp]
		\centering
		\includegraphics[width=.95\textwidth]{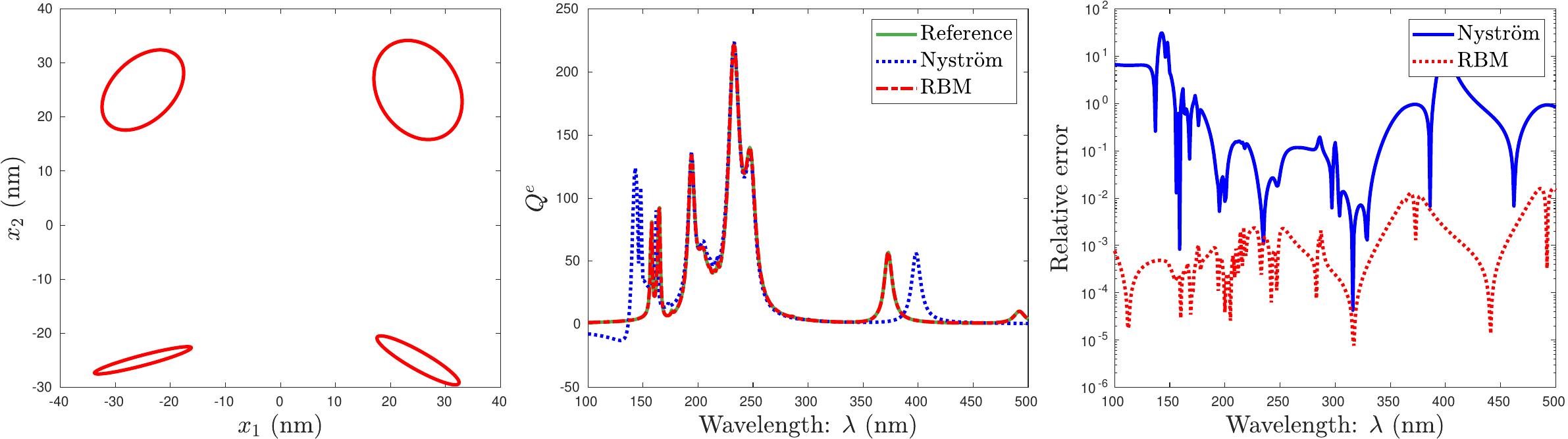}
		\caption{Left: Configuration of multiple particles ($M = 4$) with random parameters. Center: Extinction cross section $Q^e$. Right: Relative error $|Q^e - Q^e_{\mathrm{ref}}| / |Q^e_{\mathrm{ref}}|$. The number of eigenfuncions $N = 10$ for each particle.}
		\label{fig:scattering_multiple}
	\end{figure}

	Next, we consider a configuration involving multiple randomly distributed particles, as illustrated in \cref{fig:scattering_multiple}. In this case, RBM continues to outperform the Nystr\"{o}m method, as shown in \cref{fig:scattering_multiple}. RBM accurately captures the reference solution and maintains low relative errors using only a few eigenfunctions, thereby substantially reducing computational cost.

\begin{figure}[htbp]
		\centering
		\includegraphics[width=.95\textwidth]{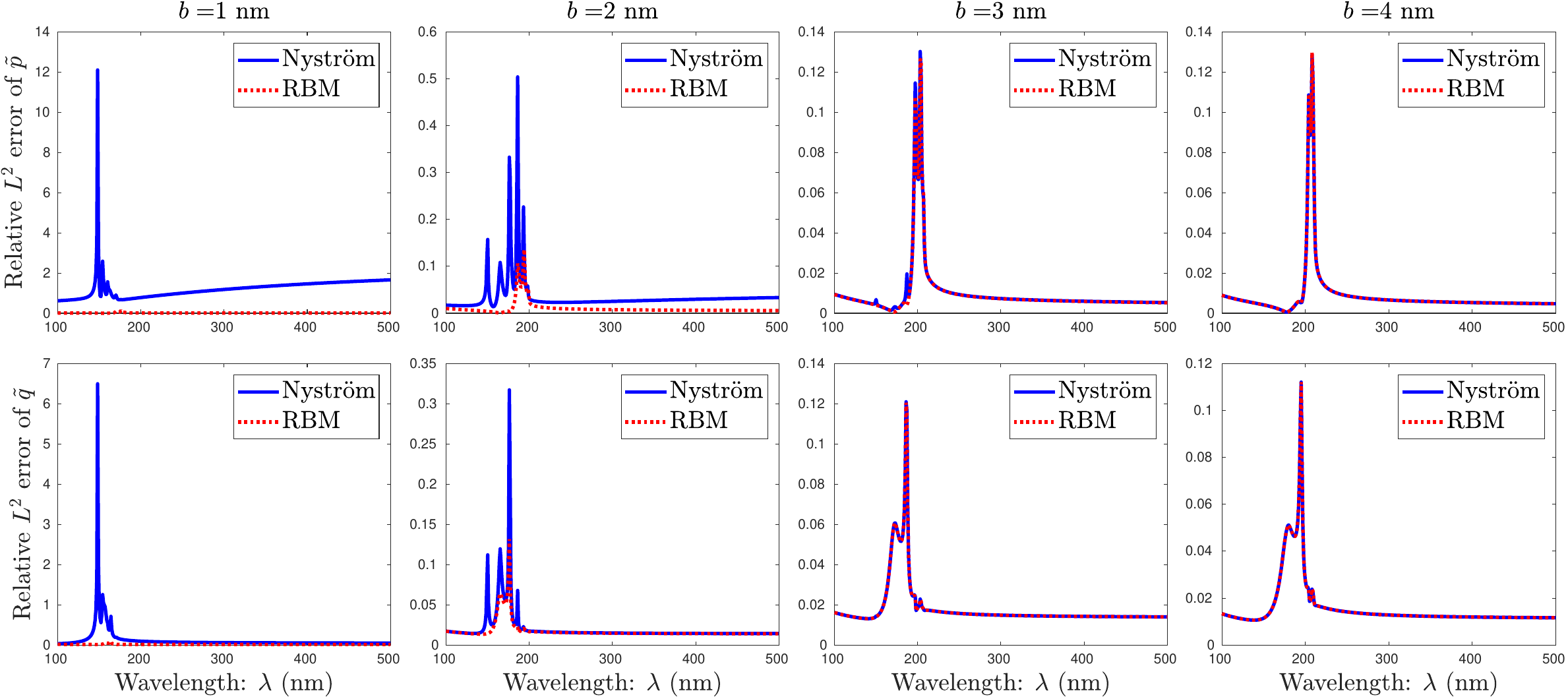}
		\caption{Relative $L^2$ errors of $\tilde{p}$ and $\tilde{q}$ for a single particle located at the origin with $a = 10$~nm, $\theta = 0$, and the number of eigenfuncions $N = 10$. }
		\label{fig:adjoint_single}
	\end{figure}

\begin{figure}[htbp]
		\centering
		\includegraphics[width=.95\textwidth]{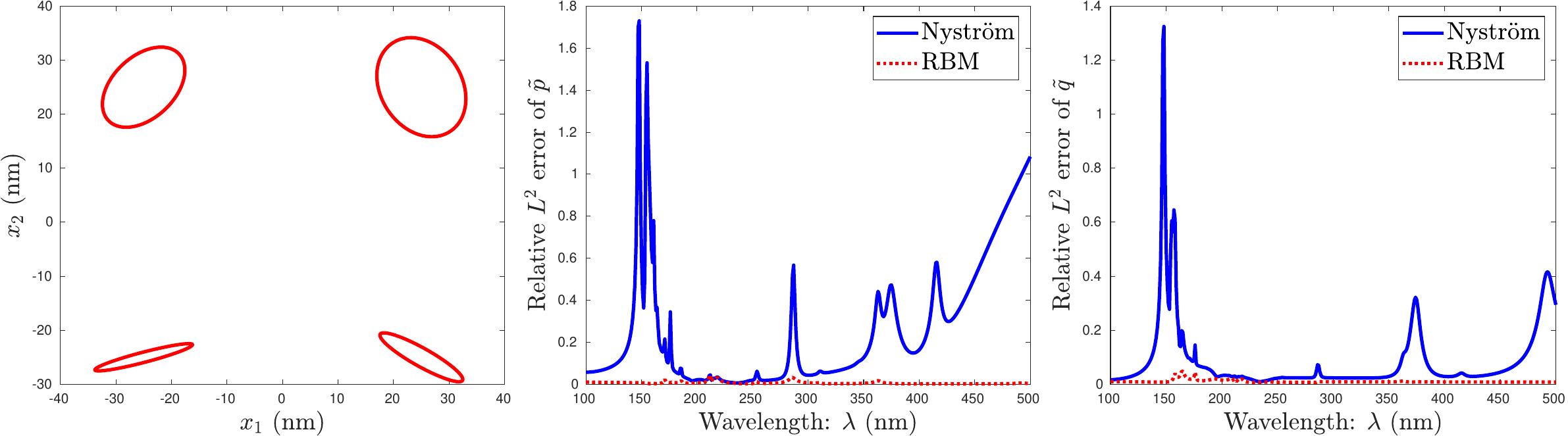}
		\caption{Left: Configuration of multiple particles ($M = 4$) with random parameters. Center: Relative $L^2$ error of $\tilde{p}$. Right: Relative $L^2$ error of $\tilde{q}$. The number of eigenfuncions $N = 10$ for each particle. }
		\label{fig:adjoint_multiple}
	\end{figure}

	Second, we evaluate the efficiency of RBM for solving the adjoint problem in comparison with the Nystr\"{o}m method. For both the single-particle case (\cref{fig:adjoint_single}) and the multiple-particle case (\cref{fig:adjoint_multiple}), RBM outperforms the Nystr\"{o}m method, particularly when $a \gg b$, where the Nystr\"{o}m method exhibits very large relative errors.
	
\subsection{Initial guess generation}\label{example:initial_guess}
	
\begin{figure}[htbp]
	\centering
	\includegraphics[width=.95\textwidth]{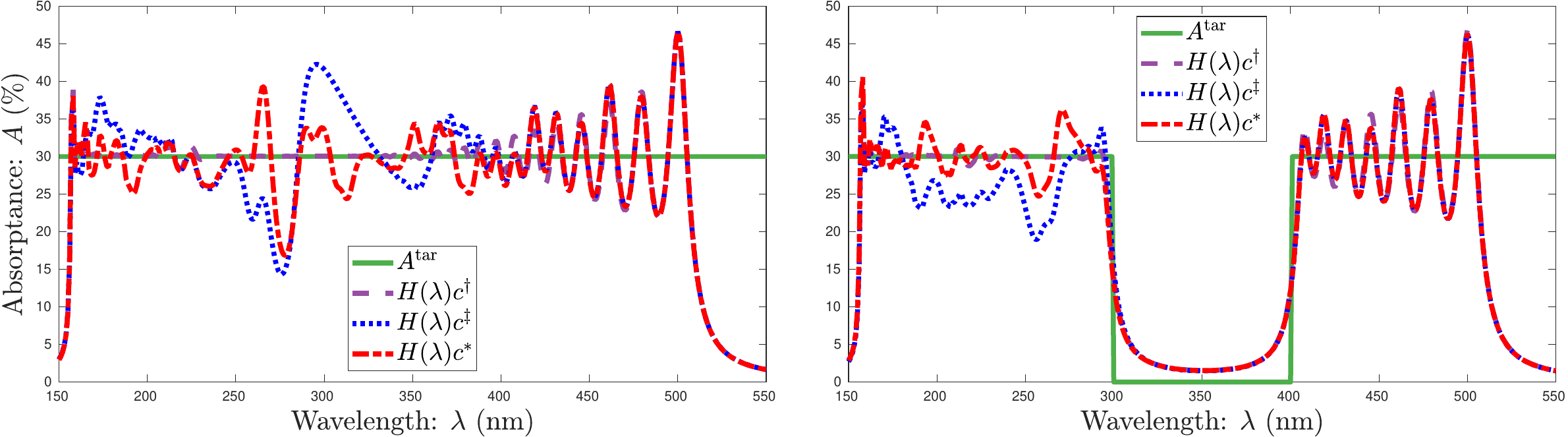}
	\caption{Initial guess generated by \cref{algo:initial_guess}. Left: constant absorptance with total number of particles $M = 104$; Right: non-constant absorptance with total number of particles $M = 83$.}
	\label{fig:initial}
\end{figure}

	In this example, we demonstrate the implementation of \cref{algo:initial_guess}. The dataset is generated by fixing $a = 10$~nm and uniformly discretizing the interval $b \in [1, 9]$~nm into 80 grid points and $\theta \in [0, \pi/2]$ into 40 grid points, resulting in a total number of  samples $L = 2511$. We then compute the absorptance $A(w_\ell, \lambda)$ over the wavelength range $\lambda \in [150, 550]$~nm. In the following, we consider both  constant and non-constant target absorptance cases, defined by
	\begin{itemize}[leftmargin=*]
		\item Constant target absorptance: $A^{\mathrm{tar}}(\lambda) = 30\%$ for $\lambda \in [150, 550]$~nm.
		\item Non-constant target absorptance: $A^{\mathrm{tar}}(\lambda) = 30\%$ for $\lambda \in [150, 300] \cup [400, 550]$~nm, and $A^{\mathrm{tar}}(\lambda) = 0$ for $\lambda \in [300, 400]$~nm.
	\end{itemize}
	\Cref{fig:initial} shows the numerical results for the constant and non-constant target absorptance cases, respectively. The relaxed solution $H(\lambda)c^\dagger$ fits the target absorptance $A^{\mathrm{tar}}$ very well using quadratic programming. The rounded solution $H(\lambda)c^\ddagger$ exhibits greater oscillation compared to $H(\lambda)c^\dagger$. Next, we employ a heuristic algorithm—particle swarm optimization—using $c^\ddagger$ as the initial guess to obtain $c^*$. As shown, $H(\lambda)c^*$ provides a better approximation to the target absorptance than $H(\lambda)c^\ddagger$.

\subsection{Optimal design}\label{example:optimal_design}
	
	In this example, we present numerical results for the optimal design problem solved using \cref{algo:optimal_design} with the following parameters:
	\begin{itemize}[leftmargin=*]
		\item Number of iterations: $N_{\mathrm{iter}} = 1000$;
		\item Number of wavelength discretization points: $N_{\lambda} = 400$;
		\item Semi-major axis range: $a_{\min} = 8$~nm, $a_{\max} = 20$~nm;
		\item Step size: $\beta = 0.2$ and aspect ratio bounds: $\eta_{\min} = 0.1$, $\eta_{\max} = 0.9$.
	\end{itemize}
	
\begin{figure}[htbp]
		\centering
		\includegraphics[width=.95\textwidth]{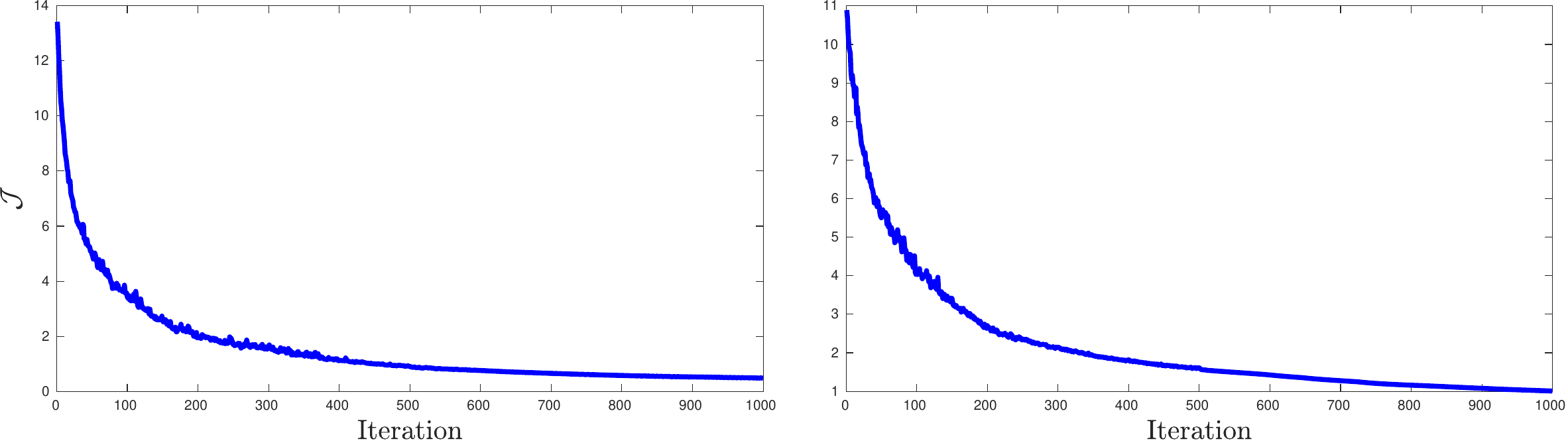}
		\caption{Objective function versus iteration for the optimal design. Left: constant target absorptance; right: non-constant target absorptance.}
		\label{fig:optimal_loss}
	\end{figure}

\begin{figure}[htbp]
		\centering
		\includegraphics[width=.95\textwidth]{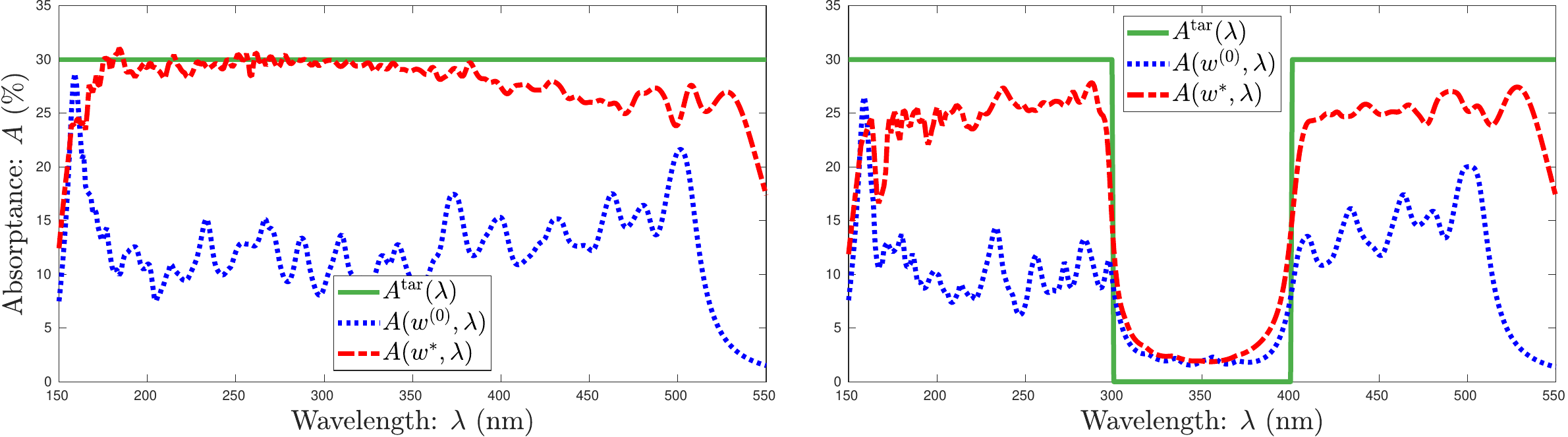}
		\caption{Comparison of target, initial, and optimized absorptance for the optimal design. Left: constant target absorptance; right: non-constant target absorptance.}
		\label{fig:optimal_absorptance}
	\end{figure}

\begin{figure}[htbp]
		\centering
		\includegraphics[width=.495\textwidth]{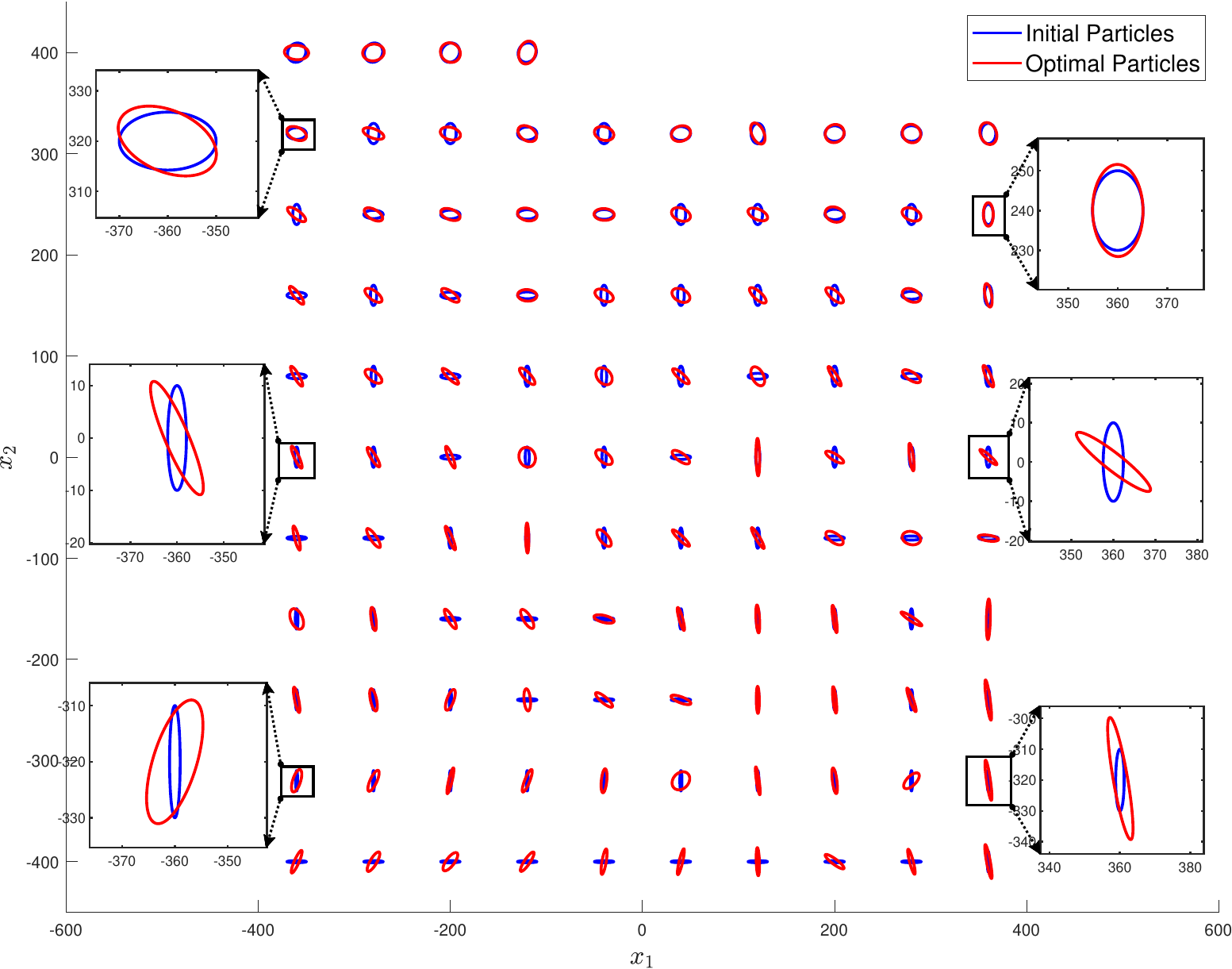}
		\includegraphics[width=.495\textwidth]{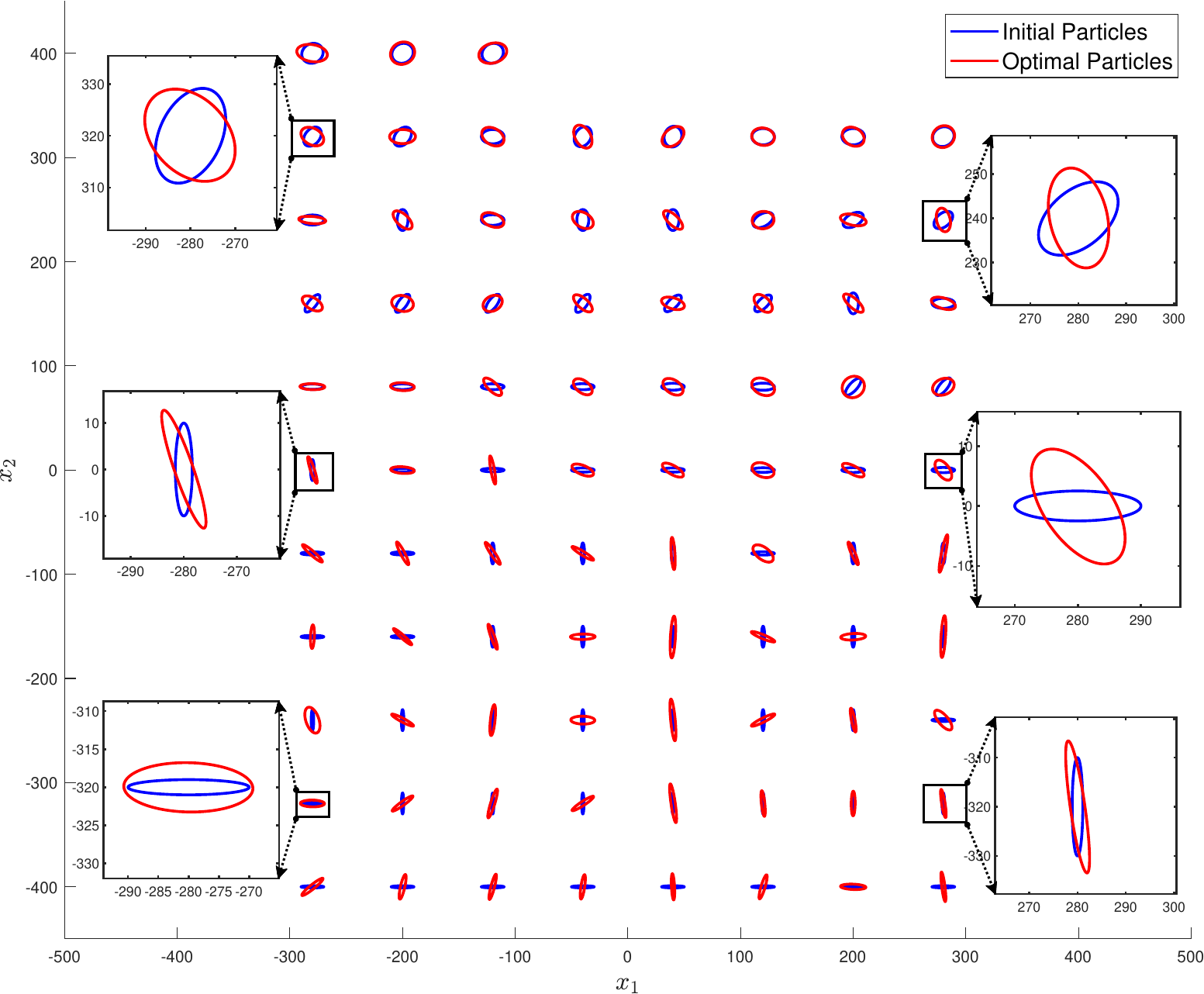}
		\caption{Initial and optimized shapes of multiple particles. Left: constant target absorptance (number of particles $M = 104$); right: non-constant target absorptance (number of particles $M = 83$).}
		\label{fig:optimal_particles_constant}
	\end{figure}

	First, \cref{fig:optimal_loss} shows that the objective function decreases and plateaus as the algorithm converges. Second, as shown in \cref{fig:optimal_absorptance}, the optimized nanoparticles achieve absorptance values very close to the target and show substantial improvement over the initial guess (see \cref{example:initial_guess}). Third, the corresponding optimized particle shapes are illustrated in \cref{fig:optimal_particles_constant}.
	
\section{Conclusions}\label{sec:conclusion}

This work presents an optimal design framework that integrates the RBM for solving problems involving multiple plasmonic nanoparticles with a physics-informed initial guess for gradient descent optimization. The optimal design problem is reformulated in a parameterized form by introducing shape parameters for multiple particles, and the corresponding shape derivatives are derived using the adjoint method. The scattering and adjoint problems are efficiently solved using distinct eigenfunctions of the NP operator and its adjoint. Notably, these eigenfunctions are shape-dependent and facilitate the efficient handling of singularities through operator splitting. Furthermore, the physics-informed initial guess is constructed under the weak scattering assumption, enhancing the convergence and accuracy of the optimization process.

The proposed optimization framework opens several promising directions for future research in material design, particularly in the context of periodic metamaterials and three-dimensional scattering problems governed by Maxwell’s equations. To further enhance broadband absorption, it is crucial to incorporate more geometrically complex resonant structures—such as nearly touching particles \cite{yu2018plasmonic} and crescent-shaped particles \cite{aubry2010broadband}, as resonance behavior is highly sensitive to particle shape. However, the analysis and simulation of such intricate geometries will necessitate the development of new, efficient numerical methods. Moreover, future design efforts should aim to optimize trade-offs between absorption performance and the total volume of the particles. Achieving this balance requires a deeper understanding of the fundamental relationship between broadband absorption and particle volume, constrained by physical limits inherent to passive systems; see, for example, \cite{gustafsson2011physical,meng2022fundamental,nedic2021herglotz,yang2017optimal}.

\appendix
	
\section{Proof of Theorem~\ref{thm:asymptotic_of_energy_flow}}\label{appendix:proof_of_energy_flow}

	Using the single-layer potential representation from \eqref{eq:integral_representation}, the scattered field $u^s$ admits the asymptotic expansion:
	\begin{align*}
		u^s(x) = \mathcal{S}^{k_m}[\varphi](x) = \int_{\partial D} G^{k_m}(x - y) \varphi(y) \, \mathrm{d}\sigma(y)
		= \frac{e^{\mathrm{i} k_m |x|}}{\sqrt{|x|}} \left\{ u^{\infty}(\hat{x}) + \mathcal{O}\left(\frac{1}{|x|}\right) \right\},
	\end{align*}
	where $u^{\infty}(\hat{x})$ is the far-field pattern defined in \eqref{eq:far_field_operator}. The corresponding gradient field exhibits the following asymptotic behavior:
	\begin{align*}
		\nabla u^s(x) = \frac{e^{\mathrm{i} k_m |x|}}{\sqrt{|x|}} \left\{ \mathrm{i} k_m \hat{x} \, u^{\infty}(\hat{x}) + \mathcal{O}\left(\frac{1}{|x|}\right) \right\}.
	\end{align*}
	Based on the energy flux definitions in \eqref{eq:def_of_Flux}, we derive the following asymptotic expansions through direct calculations:
	\begin{align*}
		F^i &= 2C k_m d, \quad
		F^s = 2C k_m \frac{|u^{\infty}|^2}{|x|} \hat{x} + \mathcal{O}\left(\frac{1}{|x|^2}\right), \\
		F^\prime &= 2C k_m \Im\left( \frac{e^{\mathrm{i} k_m |x| (1 - d \cdot \hat{x})}}{\sqrt{|x|}} \, \mathrm{i} (d + \hat{x}) \, u^{\infty}(\hat{x}) \right) + \mathcal{O}\left(\frac{1}{|x|}\right).
	\end{align*}
	Applying the definition of energy flow from \eqref{eq:def_of_Flow}, we derive the following asymptotic expansions:
    \begin{equation*}
    \resizebox{\hsize}{!}{$
        	\begin{aligned}
		E^i_{R,\Theta}
		&= \int_{\partial B_{R,\Theta}} 2C k_m d \cdot \nu(x) \, \mathrm{d}\sigma(x) = 2C k_m \int_{\bar{\theta} - \Delta\theta}^{\bar{\theta} + \Delta\theta} \cos(\theta - \theta_0) \, \mathrm{d}\theta = 4C k_m R \sin(\Delta\theta) \cos(\bar{\theta} - \theta_0), \\
		E^s_{R,\Theta}
		&= \int_{\partial B_{R,\Theta}} \left[ 2C k_m \frac{|u^{\infty}|^2}{R} \hat{x} + \mathcal{O}\left(\frac{1}{R^2}\right) \right] \cdot \nu(x) \, \mathrm{d}\sigma(x) = 2C k_m \int_{\bar{\theta} - \Delta\theta}^{\bar{\theta} + \Delta\theta} \left| u^{\infty}(\hat{x}(\theta)) \right|^2 \, \mathrm{d}\theta + \mathcal{O}\left(\frac{1}{R}\right), \\
		E^\prime_{R,\Theta}
		&= \int_{\partial B_{R,\Theta}} \left[ 2C k_m \Im\left( \frac{e^{\mathrm{i} k_m R (1 - d \cdot \hat{x})}}{\sqrt{R}} \, \mathrm{i} (d + \hat{x}) \, u^{\infty}(\hat{x}) \right) + \mathcal{O}\left(\frac{1}{R}\right) \right] \cdot \nu(x) \, \mathrm{d}\sigma(x) \\
		&= 2C k_m \sqrt{R} \Im\left[ \int_{\bar{\theta} - \Delta\theta}^{\bar{\theta} + \Delta\theta} e^{\mathrm{i} k_m R (1 - d \cdot \hat{x}(\theta))} \, \mathrm{i} (d + \hat{x}(\theta)) \, u^{\infty}(\hat{x}(\theta)) \, \mathrm{d}\theta \right] + \mathcal{O}\left(\frac{1}{\sqrt{R}}\right).
	\end{aligned}
    $}
    \end{equation*}
	To evaluate the final integral in $E^\prime_{R,\Theta}$, we denote
	\begin{equation*}
		f(\theta) = \mathrm{i} (d + \hat{x}(\theta)) \, u^{\infty}(\hat{x}(\theta)), \quad
		g(\theta) = k_m R (1 - d \cdot \hat{x}(\theta)) = k_m R (1 - \cos(\theta - \theta_0)),
	\end{equation*}
	where $d = (\cos\theta_0, \sin\theta_0)$. The phase function $g(\theta)$ has critical points at $\theta_0$ and $\theta_0 + \pi$ with properties:
	\begin{equation*}
		\begin{alignedat}{3}
			f(\theta_0) &= 2\mathrm{i} \, u^\infty(d), &\quad g(\theta_0) &= 0, &\quad g''(\theta_0) &= k_m R, \\
			f(\theta_0 + \pi) &= 0, &\quad g(\theta_0 + \pi) &= 2 k_m R, &\quad g''(\theta_0 + \pi) &= -k_m R.
		\end{alignedat}
	\end{equation*}
	The asymptotic behavior of $E'_{R,\Theta}$ then follows from the method of stationary phase~\cite{stein1993harmonic}.

\section{Proof of Theorem~\ref{thm:gradient}}\label{appendix:gradient}
	The directional derivative of $\mathcal{J}(w)$ with respect to $w$ in the direction $w_0$ is given by
	\begin{align*}
		\langle \mathcal{J}^\prime(w), w_0 \rangle_{W} = 2 \langle A(\tilde{\varphi}(w, \lambda), w, \lambda) - A^{\mathrm{tar}}(\lambda), dA(\tilde{\varphi}(w, \lambda), w, \lambda)(w_0) \rangle_{L^2(\Lambda)},
	\end{align*}
		where $dA(\tilde{\varphi}(w, \lambda), w, \lambda)(w_0)$ denotes the directional derivative of $A$ with respect to $w$ in the direction $w_0$. By the chain rule, this derivative can be decomposed as:
		\begin{align*}
			dA(\tilde{\varphi}(w, \lambda), w, \lambda)(w_0) = \Re \left\{ \langle A_{\tilde{\varphi}}(\tilde{\varphi}, w, \lambda), d\tilde{\varphi}(w, \lambda)(w_0) \rangle_{V} \right\} + \langle A_w(\tilde{\varphi}, w, \lambda), w_0 \rangle_{W},
		\end{align*}
		where $A_{\tilde{\varphi}}(\tilde{\varphi}, w, \lambda) \in V$ and $A_w(\tilde{\varphi}, w, \lambda) \in W$ are the partial derivatives of $A$ with respect to $\tilde{\varphi}$ and $w$, respectively. Next, we employ the adjoint method to deal with $\langle A_{\tilde{\varphi}}(\tilde{\varphi}, w, \lambda), d\tilde{\varphi}(w, \lambda)(w_0) \rangle_{V}$. By the adjoint system \eqref{eq:abstract_adjoint_equation}, we obtain:
		\begin{equation}\label{eq:dAdphi}
			\begin{aligned}
				& \langle A_{\tilde{\varphi}}( \tilde{\varphi}, w;\lambda), d\tilde{\varphi}(w,\lambda)(w_0) \rangle_{V}\\
				= &\langle A_{\tilde{\phi}}( \tilde{\varphi}, w,\lambda), d\tilde{\varphi}(w,\lambda)(w_0) \rangle_{V}+\langle A_{\tilde{\varphi}}( \tilde{\varphi}, w,\lambda), d\tilde{\varphi}(w,\lambda)(w_0) \rangle_{V}\\
				= & \langle  \tilde{p} , \mathcal{E}_{\tilde{\phi}}[\tilde{\phi}, \tilde{\varphi}, w,\lambda]d\tilde{\varphi}(w,\lambda)(w_0) +\mathcal{E}_{\tilde{\phi}}[\tilde{\phi}, \tilde{\varphi}, w, \lambda] d\tilde{\varphi}(w,\lambda)(w_0)  \rangle_{V}+ \\
				& \langle \tilde{q},  \mathcal{F}_{\tilde{\phi}}[\tilde{\phi}, \tilde{\varphi}, w, \lambda] d\tilde{\varphi}(w,\lambda)(w_0)+ \mathcal{F}_{\tilde{\varphi}}[\tilde{\phi}, \tilde{\varphi}, w,\lambda]d\tilde{\varphi}(w,\lambda)(w_0)\rangle_{V}.
			\end{aligned}
		\end{equation}
		Differentiating $\mathcal{E}$ and $\mathcal{F}$ with respect to $w$ in the direction $w_0$, we obtain:
		\begin{equation*}\label{eq:dEdF}
			\begin{cases}
				\mathcal{E}_{\tilde{\phi}}[\tilde{\phi}, \tilde{\varphi}, w, \lambda] d\tilde{\varphi}(w, \lambda)(w_0) + \mathcal{E}_{\tilde{\varphi}}[\tilde{\phi}, \tilde{\varphi}, w, \lambda] d\tilde{\varphi}(w, \lambda)(w_0) + \mathcal{E}_w[\tilde{\phi}, \tilde{\varphi}, w, \lambda] w_0 = 0, \\
				\mathcal{F}_{\tilde{\phi}}[\tilde{\phi}, \tilde{\varphi}, w, \lambda] d\tilde{\varphi}(w, \lambda)(w_0) + \mathcal{F}_{\tilde{\varphi}}[\tilde{\phi}, \tilde{\varphi}, w, \lambda] d\tilde{\varphi}(w, \lambda)(w_0) + \mathcal{F}_w[\tilde{\phi}, \tilde{\varphi}, w, \lambda] w_0 = 0.
			\end{cases}
		\end{equation*}
		Substituting the above equations into \eqref{eq:dAdphi}, we obtain:
		\begin{equation*}
			\begin{aligned}	
				&\langle A_{\tilde{\varphi}}(\tilde{\varphi}, w, \lambda), d\tilde{\varphi}(w, \lambda)(w_0) \rangle_{V} = -\langle \tilde{p}, \mathcal{E}_w[\tilde{\phi}, \tilde{\varphi}, w, \lambda] w_0 \rangle_{V} - \langle \tilde{q}, \mathcal{F}_w[\tilde{\phi}, \tilde{\varphi}, w, \lambda] w_0 \rangle_{V}, \\
				= & -\langle \mathcal{E}^*_w[\tilde{\phi}, \tilde{\varphi}, w, \lambda] \tilde{p} + \mathcal{F}^*_w[\tilde{\phi}, \tilde{\varphi}, w, \lambda] \tilde{q}, w_0 \rangle_{W}.
			\end{aligned}
		\end{equation*}
		Thus, we have the expression of $\langle \mathcal{J}^\prime(w),w_0 \rangle_{W}$ in the following
		\begin{equation*}
			\begin{aligned}
				= & 2\langle A(\tilde{\varphi},w, \lambda)-A^{\mathrm{tar}}(\lambda), \langle A_w( \tilde{\varphi}, w,\lambda),w_0 \rangle_{W}-\Re\left(\langle \mathcal{E}^*_w[\tilde{\phi},\tilde{\varphi}, w, \lambda]\tilde{p}+\mathcal{F}^*_w[\tilde{\phi},\tilde{\varphi}, w, \lambda]\tilde{q} , w_0  \rangle_{W}\right) \rangle_{L^2(\Lambda)},\\
				= & \langle 2\langle A(\tilde{\varphi},w, \lambda)-A^{\mathrm{tar}}(\lambda),  A_w( \tilde{\varphi}, w,\lambda) -
				\Re\left( \mathcal{E}^*_w[\tilde{\phi},\tilde{\varphi}, w, \lambda]\tilde{p}+\mathcal{F}^*_w[\tilde{\phi},\tilde{\varphi}, w, \lambda]\tilde{q} \right) \rangle_{L^2(\Lambda)},w_0 \rangle_{W} ,\\
			\end{aligned}
		\end{equation*}
		where the final equality is obtained by changing the order of integration. Therefore, 
		\begin{equation*}
			\mathcal{J}^\prime(w) = 2 \left\langle A(\tilde{\varphi}, w, \lambda) - A^{\mathrm{tar}}(\lambda), \ A_w(\tilde{\varphi}, w, \lambda) - \Re \left( \mathcal{E}^*_w[\tilde{\phi}, \tilde{\varphi}, w, \lambda] \tilde{p} + \mathcal{F}^*_w[\tilde{\phi}, \tilde{\varphi}, w, \lambda] \tilde{q} \right) \right\rangle_{L^2(\Lambda)}.
		\end{equation*}
		The proof is complete by rewriting the result in the inner product form.
		
	\section{Proof of Theorem~\ref{thm:spectral_of_Kwm}}\label{appendix:spectral_of_Kwm}
	
	Let $D$ be an ellipse centered at the origin, parameterized by $w$ with foci at $c$. According to the results in \cite{AK2016Analysis,chung2014cloaking}, the Green's function for Laplace's equation can be expressed as an eigenfunction expansion utilizing the properties of the NP operator:
	\begin{align*}\label{eq:expansion_of_G}
		G(x-y)=\frac{\rho_x+\ln(c/2)}{2\pi}-\sum_{n=1}^{\infty}\frac{\cos(nt)\cosh(n\rho)\cos(ns)+\sin(nt)\sinh(n\rho)\sin(ns) }{\pi ne^{n\rho_x}},
	\end{align*}
	where $x = (\rho_x, t) \in \mathbb{R}^2 \setminus \overline{D}$ and $y = (\rho, s) \in \partial D$ are represented in elliptic coordinates as defined in \cref{eq:elliptic_coordinate}. Utilizing the expansion of the Green's function, we obtain:
	\begin{align*}
		S^*_w[\cos(is)](t) = -\frac{\cosh(i\rho)\cos(nt)}{ie^{i\rho_x}}\Xi(\rho_x,t),\
		S^*_w[\sin(is)](t) = -\frac{\sinh(i\rho)\sin(nt)}{ie^{i\rho_x}}\Xi(\rho_x,t),
	\end{align*}
	where we have used the fact that $|x'(t)| = \Xi(\rho_x,t)$. In the special case when $i = 0$, $S^*_w[\cos(is)](t) = (\rho_x+\ln(c/2))\Xi(\rho_x,t)$.
	Note that $S_w^*[\psi]$ is continuous as $x \to \partial D$ due to the weak singularity of $G(x - y)$. Taking the limit $x \to \partial D$ (i.e., $\rho_x \to \rho$), then
	\begin{align*}\label{eq:Sswpsi}
		S^*_w[\cos(is)](t) = -\frac{\cosh(i\rho)\cos(it)}{ie^{i\rho}}\Xi(\rho,t),\
		S^*_w[\sin(is)](t) = -\frac{\sinh(i\rho)\sin(it)}{ie^{i\rho}}\Xi(\rho,t).
	\end{align*}
	
	According to the normal derivative in elliptic coordinates \cref{eq:differential_property_of_ellipse}, we have:
	\begin{equation}\label{eq:dGdnuy_in_elliptic_coordinate}
		\frac{\partial G(x-y)}{\partial \nu_y} = -\sum_{n=1}^{\infty} \frac{\cos(nt)\sinh(n\rho)\cos(ns) + \sin(nt)\cosh(n\rho)\sin(ns)}{\pi \Xi(\rho,s) e^{n\rho_x}}.
	\end{equation}
	Let $\mathcal{D}_w$ denote the adjoint parameterized form of the double layer potential:
	\begin{equation*}
		\mathcal{D}_w[\psi](t)= \int_{0}^{2\pi} \frac{\partial G(x(t)-y(s))}{\partial \nu_y} |x'(t)| \psi(s) \, \mathrm{d}s, \ x = (\rho_x, t) \in \mathbb{R}^2 \setminus \overline{D}, \ y = (\rho, s) \in \partial D.
	\end{equation*}
	Using the expansion in \cref{eq:dGdnuy_in_elliptic_coordinate}, we directly obtain:
	\begin{align*}
		\mathcal{D}_w[q_{i}^s](t) = -\frac{\cosh(i\rho)\sin(it)}{e^{i\rho_x}} \Xi(\rho_x,t),\
		\mathcal{D}_w[q_{i}^c](t) = -\frac{\sinh(i\rho)\cos(it)}{e^{i\rho_x}} \Xi(\rho_x,t),
	\end{align*}
	where $q_{i}^s = \sin(it)\Xi(\rho,t)$ and $q_{i}^c = \cos(it)\Xi(\rho,t)$. Consequently, applying the jump formula $\mathcal{D}_w[\psi]\big{|}_+ = \left(-\frac{1}{2}\mathcal{I} + \mathcal{K}_w\right)[\psi]$, we derive the explicit expressions for $\mathcal{K}_w$:
	\begin{equation*}\label{eq:Kwpsi}
		\begin{aligned}
			\mathcal{K}_w[q_{i}^s](t) &= -\frac{\cosh(i\rho)\sin(it)}{e^{i\rho}} \Xi(\rho,t) + \frac{1}{2}\sin(it)\Xi(\rho,t) = \frac{-1}{2e^{2i\rho}}q_{i}^s(t),\\
			\mathcal{K}_w[q_{i}^c](t) &= -\frac{\sinh(i\rho)\cos(it)}{e^{i\rho}} \Xi(\rho,t) + \frac{1}{2}\cos(it)\Xi(\rho,t) = \frac{1}{2e^{2i\rho}}q_{i}^c(t).
		\end{aligned}
	\end{equation*}
	
	\section{Singular behavior of Green's Function}\label{appendix:asymptotic_of_Green}

	First, we consider the asymptotic form of some radial function. According to the defnition of Bessel functions in \cite{CK2019Inverse}, we have the following asymptotic expansion
		\begin{align}
			-\frac{\mathrm{i}}{4} H_0^{(1)}(kr) & 
			=\frac{1}{2 \pi} \ln r-\left(\frac{\mathrm{i}}{4}-\frac{\gamma}{2 \pi}-\frac{1}{2 \pi} \ln \frac{k}{2}\right)+O(r^2\ln r),\label{eq:ap_gkr}\\
			\frac{\mathrm{i}k}{4} H_1^{(1)}(kr)  
			& = 	\frac{1}{2 \pi r}-\frac{k^2}{4\pi}r\ln r+\frac{k^2}{2}(\frac{\mathrm{i}}{4}-\frac{\gamma}{2 \pi}-\frac{1}{2 \pi} \ln \frac{k}{2}+\frac{1}{4\pi})r +O(r^2\ln r).\label{eq:ap_dgkrdr}
		\end{align}
	For the convenience of subsequent discussions, we denote
	\begin{equation}\label{eq:ap_radius_green_fun}
		g(r) := \frac{1}{2\pi}\ln r, \quad
		g(r;k) := -\frac{\mathrm{i}}{4} H_0^{(1)}(kr), \quad
		\widehat{g}(r; k) := g(r;k) - g(r).
	\end{equation}
	Then, we have $G(x-y)=g(|x-y|)$, $G^k(x-y)=g(|x-y|;k)$, and $\widehat G^k(x-y)=\widehat{g}(|x-y|;k)$. Let $\widehat{g}'(r; k)$ denotes the derivative with respect to $r$, using the asymptotic expansions of the Hankel function from \cref{eq:ap_gkr,eq:ap_dgkrdr}, we obtain:
	\begin{equation}\label{eq:asymptotic_of_g}
		\begin{aligned}
			\widehat{g}(r; k) &= -\left(\frac{\mathrm{i}}{4} - \frac{\gamma}{2\pi} - \frac{1}{2\pi}\ln\frac{k}{2}\right) + O(r^2\ln r), \\
			\widehat{g}'(r; k) &= -\frac{k^2}{4\pi}r\ln r + \frac{k^2}{2}\left(\frac{\mathrm{i}}{4} - \frac{\gamma}{2\pi} - \frac{1}{2\pi}\ln\frac{k}{2} + \frac{1}{4\pi}\right)r+ O(r^2\ln r),
		\end{aligned}
	\end{equation}
	which implies that $\widehat{g}(0; k) $ and $\widehat{g}^\prime(r; k)$ is non-singular at $0$. Thus, we can define the limiting values: $	\widehat{g}(0; k):= -\left(\frac{\mathrm{i}}{4} - \frac{\gamma}{2\pi} - \frac{1}{2\pi}\ln\frac{k}{2}\right)$ and $\widehat{g}^\prime(0; k):= 0$. Notice the derivative relation for Hankel functions:
	\begin{align*}
		H_1^{\prime(1)}(z) = H_0^{(1)}(z) - \frac{1}{z}H_1^{(1)}(z) \implies H_1^{\prime(1)}(kr) = kH_0^{(1)}(kr) - \frac{1}{r}H_1^{(1)}(kr).
	\end{align*}
	This leads to the asymptotic behavior of the second derivative:
		\begin{align*}\label{eq:appendix_hatgpp}
		\widehat{g}^{\prime\prime}(r; k) = -k^2g(r;k)-\frac{\widehat{g}^{\prime}(r; k)}{r}= -\frac{k^2}{4\pi}\ln r+\frac{k^2}{2}(\frac{\mathrm{i}}{4}-\frac{\gamma}{2 \pi}-\frac{1}{2 \pi} \ln \frac{k}{2}-\frac{1}{4\pi})+O(r\ln r),
		\end{align*}
	which shows that $\widehat{g}^{\prime\prime}(r; k)$ is singular at $0$.
	
	Before going further, we give the properties of parameterized elliptic curve. By direct calculation, we have the following results
	\begin{equation}\label{eq:property_of_hatz}
		\begin{aligned}
			\widehat{z}(t,s;w_m)\cdot \frac{\partial z(t,s;w_m)}{\partial w} & \sim |t-s| ( a_m\sin^2t,  b_m\cos^2t, 0,  0,  0)/|x^\prime(t;w_m)|,\\
			\widehat{z}(t,s;w_m)\cdot \nu_x(t;w_m) & \sim |t-s|\frac{a_mb_m}{2|x^{\prime}(t;w_m)|^2},\\
			\frac{\partial (\widehat{z}(t,s;w_m)\cdot \nu_x(t;w_m))}{\partial w} &\sim |t-s|\frac{a_m^4\sin^4t-b_m^4\cos^4 t}{2|x^{\prime}(t;w_m)|^6}( -b_m, a_m, 0,  0,  0),
		\end{aligned}
	\end{equation}
	where $\nu_x(t;w_m)= -\mathrm{i}x^{\prime}(t;w_m)/|x^\prime(t;w_m)|$ and $\widehat{z}(t,s;w_m) =  z(t,s;w_m)/|z(t,s;w_m)|$.

	Based on the preceding preparations, we now present the singular properties of the Green's function. According to \cref{eq:ap_radius_green_fun,eq:asymptotic_of_g}, it is straightforward to obtain
	\begin{equation}\label{eq:limit_of_hatGk}
		\begin{aligned}
			\lim\limits_{t\rightarrow s}\widehat G^{k}(z(t,s;w_m)) = & \lim\limits_{t\rightarrow s} \widehat{g}(|z(t,s;w_m)|; k)=\widehat{g}(0; k),\\
			\lim\limits_{t\rightarrow s} \frac{\partial G^k(z(t,s;w_m))}{\partial \nu_x(t;w_m)} = & \lim\limits_{t\rightarrow s} \widehat{g}^{\prime}(|z(t,s;w_m)|;k) \widehat{z}(t,s;w_m)\cdot \nu_x(t;w_m) = 0.
		\end{aligned}
	\end{equation}
	Similar argument, we also have
	\begin{align*}
			\lim\limits_{t\rightarrow s}\widehat G^{-k}(z(t,s;w_m)) =  \widehat{g}(0; -k)\quad \text{and}	
			\quad \lim\limits_{t\rightarrow s} \frac{\partial \widehat G^{-k}(z(t,s;w_m))}{\partial \nu_y(s;w_m)}  = 0.
	\end{align*}

	Next, we consider the derivative of $\widehat G^{k}(z(t,s;w_m))$ and $\frac{\partial G^k(z(t,s;w_m))}{\partial \nu_x(t;w_m)}$ with respect $w_m$ term by term. By using the property  \eqref{eq:property_of_hatz}, we have
	\begin{equation*}\label{eq:limit_of_dhatGkdw}
		\begin{aligned}
			\lim\limits_{t\rightarrow s}\frac{\partial (\widehat G^{k}(z(t,s;w_m)))}{\partial w_m} 
			= \lim\limits_{t\rightarrow s} \widehat{g}^{\prime}(|z(t,s;w_m)|;k)\widehat{z}(t,s;w_m)\cdot \frac{\partial z(t,s;w_m)}{\partial w}=( 0,  0, 0,  0,  0)^\top,\\
		\end{aligned}
	\end{equation*}
	\begin{equation*}\label{eq:limit_of_dhatGkdnudw}
		\begin{aligned}
			&\lim\limits_{t\rightarrow s} \frac{\partial}{\partial w_m} \left(\frac{\partial \widehat{G}^{k}(z(t,s;w_m))}{\partial \nu_x(t;w_m)} \right)= \lim\limits_{t\rightarrow s} \frac{\partial}{\partial w_m} \left( \widehat{g}^{\prime}(|z(t,s;w_m)|;k)\widehat{z}(t,s;w_m)\cdot \nu_x(t;w_m) \right)\\
			=   & \lim\limits_{t\rightarrow s} \left[ \widehat{g}^{\prime\prime}(|z(t,s;w_m)|;k)\widehat{z}(t,s;w_m)\cdot \frac{\partial z(t,s;w_m)}{\partial w_m} \widehat{z}(t,s;w_m)\cdot \nu_x(t;w_m)+\right.\\
			&\quad \quad \left.\widehat{g}^{\prime}(|z(t,s;w_m)|;k)\frac{\partial (\widehat{z}(t,s;w_m)\cdot \nu_x(t;w_m))}{\partial w_m}\right]
			=( 0,  0, 0,  0,  0)^\top.
		\end{aligned}
	\end{equation*}

	\section{Derivative of Boundary Operators}\label{sec:appendix_Derivative_of_BO}
	
	For convenience in calculations, let $w = (a, b, \theta, x_1, x_2)$ be the parameters of a single ellipse (if multiple ellipses are considered, we use subscript $m$, i.e., $w_m$ for the $m$-th ellipse). According to the \cref{thm:spectral_of_NP}, we have
	\begin{equation}\label{eq:dpsidwm}
		\begin{aligned}
			\frac{\partial \psi_{i}}{\partial w}  = -\frac{1}{\Xi(\rho,t)}\frac{\partial \Xi(\rho,t)}{\partial w}\psi_{i} = -(\frac{1}{c}\frac{\partial c}{\partial w}+\frac{ab}{\Xi^2(\rho,t)}\frac{\partial \rho}{\partial w})\psi_{i}.
		\end{aligned}
	\end{equation}
	where $\psi_{i}$ denotes either $\psi^s_{i}$ or $\psi^c_{i}$ and
	\begin{equation*}
		\begin{aligned}
			\frac{\partial c}{\partial w}  = (0, a, -b, 0, 0)^\top/c,\quad \frac{\partial \rho}{\partial w}  = (0, -b, a, 0, 0)^\top/c^2,\quad \frac{\partial \alpha_{i}}{\partial w}  = -2i\alpha_{i}\frac{\partial \rho}{\partial w}.
		\end{aligned}
	\end{equation*}
	 According to the spectral property of NP operator in \cref{thm:spectral_of_NP}, we obtain:
	\begin{equation*}\label{eq:dKswpsidw}
		\begin{aligned}
			&\frac{\partial \mathcal{K}^*_{w}\left[\psi^s_{i}\right]}{\partial w} =-\frac{\partial \alpha_{i}}{\partial w} \psi^s_{i} - \alpha_{i}\frac{\partial \psi^s_{i}}{\partial w}
			= \left(\frac{1}{c}\frac{\partial c}{\partial w}+(2i+\frac{ab}{\Xi^2})\frac{\partial \rho}{\partial w}\right)\alpha_{i}\psi^s_{i},\\
			&\frac{\partial \mathcal{K}^*_{w}\left[\psi^c_{i}\right]}{\partial w} =\frac{\partial \alpha_{i}}{\partial w} \psi^c_{i} +\alpha_{i}\frac{\partial \psi^c_{i}}{\partial w}
			= -\left(\frac{1}{c}\frac{\partial c}{\partial w}+(2i+\frac{ab}{\Xi^2})\frac{\partial \rho}{\partial w}\right)\alpha_{i}\psi^c_{i},\\
			&\frac{\partial \mathcal{S}_{w}[\psi^s_{i}]}{\partial w}= \frac{\partial \alpha_{i}}{\partial w}\frac{\sin it}{i},\quad \frac{\partial \mathcal{S}_{w}[\psi^c_{i}]}{\partial w}= -\frac{\partial \alpha_{i}}{\partial w}\frac{\cos it}{i},\quad\frac{\partial \mathcal{S}_{w}[\psi_{m,0}]}{\partial w} =\frac{\partial \rho}{\partial w} +\frac{1}{c}\frac{\partial c}{\partial w}.
		\end{aligned}
	\end{equation*}

	Next, we consider the boundary operator with density $\frac{\partial \psi_{i}}{\partial w}$ in \eqref{eq:dpsidwm}:
	\begin{equation*}\label{eq:Swdpsidw_Kswdpsidw}
		\begin{aligned}
			\mathcal{S}_{w}[\frac{\partial \psi_{i}}{\partial w}] = -\frac{1}{c}\frac{\partial c}{\partial w}\mathcal{S}_{w}[\psi_{i}]-ab\frac{\partial \rho}{\partial w}\mathcal{S}_{w}[\frac{\psi_{i}}{\Xi^2}],~
			\mathcal{K}_{w}^*[\frac{\partial \psi_{i}}{\partial w}] = -\frac{1}{c}\frac{\partial c}{\partial w}\mathcal{K}_{w}^*[\psi_{i}]-ab\frac{\partial \rho}{\partial w}\mathcal{K}_{w}^*[\frac{\psi_{i}}{\Xi^2}].
		\end{aligned}
	\end{equation*}
	Note that we have explicit expressions for $\mathcal{S}_{w}[\psi_{i}]$ and $\mathcal{K}_{w}^*[\psi_{i}]$. Thus, we only need to evaluate the boundary operator with density $\psi_i/\Xi^2$. According to the expansion of Green's function in \cref{appendix:spectral_of_Kwm}, we obtain:
	\begin{equation*}
		\begin{aligned}
			\mathcal{S}_{w}[\frac{\psi^s_{i}}{\Xi^2}] &=\kappa^{cs}_{0i}\frac{\rho+\ln(c/2)}{2}-\sum_{n=1}^{\infty}\frac{\kappa^{cs}_{ni}\cos(nt)\cosh(n\rho)+\kappa^{ss}_{ni}\sin(nt)\sinh(n\rho)}{ne^{n\rho}},\\
			\mathcal{S}_{w}[\frac{\psi^c_{i}}{\Xi^2}] &=\kappa^{cc}_{0i}\frac{\rho+\ln(c/2)}{2}-\sum_{n=1}^{\infty}\frac{\kappa^{cc}_{ni}\cos(nt)\cosh(n\rho)+\kappa^{sc}_{ni}\sin(nt)\sinh(n\rho)}{ne^{n\rho}},\\
		\end{aligned}
	\end{equation*}
	where the coefficients $\kappa^{ss}_{ni}$, $\kappa^{cs}_{ni}$, $\kappa^{sc}_{ni}$, and $\kappa^{cc}_{ni}$ are defined by the following integrals:
	\begin{align*}
			\kappa^{ss}_{ni} = \int_{0}^{2\pi}  \frac{\sin(ns) \sin(is)}{c\pi(\sinh^2\rho+\sin^2s)}\mathrm{d}s,\quad
			\kappa^{cs}_{ni}  = \int_{0}^{2\pi}  \frac{\cos(ns) \sin(is)}{c\pi(\sinh^2\rho+\sin^2s)}\mathrm{d}s,\\
			\kappa^{sc}_{ni} = \int_{0}^{2\pi}  \frac{\sin(ns) \cos(is)}{c\pi(\sinh^2\rho+\sin^2s)}\mathrm{d}s,\quad
			\kappa^{cc}_{ni}  = \int_{0}^{2\pi}  \frac{\cos(ns) \cos(is)}{c\pi(\sinh^2\rho+\sin^2s)}\mathrm{d}s.
	\end{align*}
	By to the jump formula, we have $\mathcal{K}_w^*[\psi] = \frac{\partial }{\partial \nu_x}\mathcal{S}_w[\psi]|_{+}-\frac{1}{2}\psi$, which implies that
	\begin{align*}
		\mathcal{K}_w^*[\frac{\psi^s_{i}}{\Xi^2}](t) =&\frac{\kappa^{cs}_{0i}}{2}+\sum_{n=1}^{\infty}\frac{\kappa^{cs}_{ni}\cos(nt)\cosh(n\rho)+\kappa^{ss}_{ni}\sin(nt)\sinh(n\rho)}{e^{n\rho}}-\frac{1}{2}\frac{\psi^s_{i}}{\Xi^2},\\
		\mathcal{K}_w^*[\frac{\psi^c_{i}}{\Xi^2}](t) =&\frac{\kappa^{cc}_{0i}}{2}+\sum_{n=1}^{\infty}\frac{\kappa^{cc}_{ni}\cos(nt)\cosh(n\rho)+\kappa^{sc}_{ni}\sin(nt)\sinh(n\rho)}{e^{n\rho}}-\frac{1}{2}\frac{\psi^c_{i}}{\Xi^2}.
	\end{align*}

\bibliographystyle{siamplain}
\bibliography{ref}

\end{document}

%% file: ex_shared.tex
\usepackage{lipsum}
\usepackage{amsfonts}
\usepackage{graphicx}
\usepackage{epstopdf}
\usepackage{algorithmic}

\usepackage{enumitem}
\usepackage{epsfig, epstopdf}
\usepackage{tikz, tikz-cd}
\usepackage{graphicx,subfig}

\usepackage{caption}
\ifpdf
  \DeclareGraphicsExtensions{.eps,.pdf,.png,.jpg}
\else
  \DeclareGraphicsExtensions{.eps}
\fi

\newsiamremark{remark}{Remark}
\newsiamremark{hypothesis}{Hypothesis}

\crefname{hypothesis}{Hypothesis}{Hypotheses}

\newsiamthm{claim}{Claim}

\newsiamremark{fact}{Fact}
\crefname{fact}{Fact}{Facts}

\headers{Optimal Design of Multiple Nanoparticles via RBM}{Y. Gao, H. Zhang,  and K. Zhang}

\title{Optimal Design of Broadband Absorbers with Multiple Plasmonic Nanoparticles via Reduced Basis Method\thanks{\today.
		\funding{The work of Y.~Gao was partially supported by the National Key R\&D Program of China (Grant No.~2024YFA1012302). The work of H. Zhang was supported by the Hong Kong Research Grants Council General Research Fund (GRF, No.~16307024) and the National Natural Science Foundation of China (NSFC, No.~12371425). The work of K.~Zhang was partially supported by the National Natural Science Foundation of China (NSFC, No.~12271207).}}}

\author{Yu Gao\thanks{Department of Mathematics, The Hong Kong University of Science and Technology, Clear Water Bay, Hong Kong SAR, China
		(\email{yugaomath@gmail.com}).}
	\and Hai Zhang\thanks{Department of Mathematics, The Hong Kong University of Science and Technology, Clear Water Bay, Hong Kong SAR, China}
		(\email{haizhang@ust.hk})
	\and Kai Zhang\thanks{School of Mathematics, Jilin University, Changchun, China
			(\email{zhangkaimath@jlu.edu.cn})}
	}

\usepackage{amsopn}
\DeclareMathOperator{\diag}{diag}